\documentclass{amsart}
\usepackage{amsfonts,amsmath,amssymb,color,esint}
\usepackage{enumerate}
\usepackage{graphicx}
\usepackage{tikz}
\usetikzlibrary{decorations.markings}
\usetikzlibrary{plotmarks}
\usetikzlibrary{patterns}
\numberwithin{equation}{section}

\usepackage[babel=true,kerning=true]{microtype}

\usepackage{amsthm,leqno}
\usepackage{hyperref}

\newtheorem{theo}{Theorem}[section]
\newtheorem{lem}{Lemma}[section]
\newtheorem{defi}{Definition}[section]

\newtheorem{prop}{Proposition}[section]
\newtheorem{cor}{Corollary}[section]

\newcommand{\eps}{\varepsilon}
\newcommand{\D}{\mathbb{D}}
\newcommand{\R}{\mathbb{R}}

\newcommand{\ds}{\displaystyle}

\begin{document}

\title[]{Existence of Min-Max Free Boundary Disks Realising the Width of a Manifold}
\author{Paul Laurain and Romain Petrides} 
\address{Paul Laurain, Institut de Math\'{e}matiques de Jussieu, Universit\'{e} Paris Diderot, B\^{a}timent Sophie Germain\\ Case 7012, 75205 PARIS Cedex 13\\
France}
\email{paul.laurain@imj-prg.fr}
\address{Romain Petrides, Institut de Math\'{e}matiques de Jussieu, Universit\'{e} Paris Diderot, B\^{a}timent Sophie Germain\\ Case 7012, 75205 PARIS Cedex 13\\
France}
\email{romain.petrides@imj-prg.fr}

\maketitle

\begin{abstract} 
We perform a replacement procedure in order to produce a free boundary minimal surface whose area achieves the min-max value over all disk sweepouts of a manifold whose boundary lies in a submanifold. Our result is based on a proof of the convexity of the energy for free boundary harmonic maps and a generalization of Colding-Minicozzi replacement procedure.
\end{abstract}

\maketitle
\tableofcontents
\section{Introduction}

Minimal surfaces are fundamental objects in Riemannian geometry. Beyond their intrinsic interest and beauty, they have been used to prove many fundamental results in geometry and topology. Recently, the min-max approach initiated by Almgren and Pitts \cite{P81} who developed a very deep theory of existence and regularity of minimal hypersurfaces, has known a new birth in the minimal surface theory. The most striking application is the proof by Marques-Neves of the Willmore conjecture \cite{MN14} or the result by Irie-Marques-Neves who solve the Yau conjecture in the generic situation \cite{IMN18}. We have to note that Guaraco has proposed a new approach relying on the Allen–Cahn equation \cite{Gua}, see also Chodosh-Mantoulidis \cite{CM}. Rivi\`ere has also performed a promising viscosity method to solve extrinsic min-max problems for minimal surfaces in any co-dimension, see \cite{Riv17,Riv17b}.\\

In this paper we are interested in the seminal approach of the existence theory of minimal surfaces, the one used by Douglas and Rad\'o to solve the Plateau problem, namely the variational approach. In this issue, the natural continuation to look for critical points is a min-max, a setting which finds its foundation in the work of Palais \cite{Palais}. A min-max version of this approach has been performed in the very nice paper by Colding-Minicozzi which especially permits to define the {\it width} of a 3-sphere. Inspired by the work by Douglas and Rad\'o they replaced the area by the Dirichlet energy functional, and inspired by the Birkhoff's so-called curve-shortening procedure, gave a strong harmonic replacement procedure, based on the local convexity of the energy under small energy assumptions. Following the work by Colding-Minicozzi, we extend this result to free-boundary minimal surfaces. Those surfaces can be seen as a generalization of closed geodesics,  in the sense that the boundary of the disk is a $1/2$-harmonic map which is the $1$-dimensional equivalent of the harmonic map. Indeed, free boundary minimal surfaces are also subject to a non-local approach, see \cite{DP} and reference therein.\\

The main goal of this paper is to perform a free-boundary harmonic replacement procedure in order to produce a free boundary minimal surface whose area achieves the min-max value over all disk sweepouts of a manifold whose boundary lies in a submanifold.

Let $N$ be a $n$-dimensional closed manifold and $M\subset N$ a connected $m$-dimensional compact submanifold. We let
$$ \mathcal{A} = \{ u\in W^{1,2}(\mathbb{D},N)\cap\mathcal{C}^0(\overline{\mathbb{D}},N) \; \vert \; u(\partial\mathbb{D}) \subset M \} $$
be the set of admissible parametrized disks in $N$ whose boundaries lie in $M$, endowed with the $\Vert\, .\, \Vert_{L^{\infty}} + \Vert \nabla \, .\,  \Vert_{L^2}$-norm. We fix a parameter $t_0 \in \partial \mathbb{B}^{k-2}$ and an arbitrary point $m_0\in M$. We define a sweepout as a map $\sigma : \mathbb{B}^{k-2} \to \mathcal{A}$ on a $k-2$-ball which is the set of parameters with $k\geq 2$ satisfying
\begin{itemize}
\item $\forall t \in \partial \mathbb{B}^{k-2}$, $\sigma_t$ is a constant function in $M$. 
\item $t \mapsto \sigma_t$ is continuous in $\mathcal{A}$.
\item $\sigma_{t_0}=m_0$\footnote{ the condition that $\sigma_{t_0}=m_0$ can be removed when $M$ is simply connected. Else, it is necessary to get that $\pi_{k-2}(\mathcal{A})\cong \pi_k(N,M)$.}.
\end{itemize}
Let $\omega$ be a homotopy class of sweepouts, we can set a topological invariant, called the {\it width} as
\begin{equation}\label{defWarea} W(N,M,\omega)= \inf_{\sigma \in \omega} \max_{t\in \mathbb{B}^{k-2}} Area(\sigma_t). \end{equation}
Notice that there is a non-trivial homotopy class of sweepouts (or there is a sweepout non-homotopic to a constant one) if and only if $\pi_k(N,M) \neq \{0\}$, see section 1.3 of \cite{F00}. In particular in those cases, $W(N, M, \omega)>0$ and our main theorem applies. Our setting is very general since it contains the classical Plateau problem for $M$ being a closed curve and $k=2$. And it of course contains also some real min-max, for instance, if $M$ is the boundary of a strictly-convex domain in $\mathbb{R}^3$, the level set of the height function generates a non-trivial homotopy class and $\pi_3(N,M) \neq \{0\}$.\\

The definition of the classical {\it width} goes back to Birkhoff \cite{B} for the problem of geodesics : it is the smallest length we need for some circle to pull-over a compact manifold. In our context, we define the smallest area we need to cross over a compact manifold with an interface which has the topology of a disk whose boundary slides on this compact manifold.\\

Like for the problem of geodesics or the Douglas-Rad\'o approach, it is more convenient to use the energy as a functional instead of the area. We know that for $u\in \mathcal{A}$, $Area(u)\leq E(u):= \frac{1}{2}\int_{\mathbb{D}} \vert\nabla u\vert^2$ with equality if and only if $u : \mathbb{D}\to N$ is conformal. Moreover, on a disk, for any immersion $u\in \mathcal{A}$, we can change the parametrization so that $u \circ \phi \in \mathcal{A}$ is a "almost-conformal" map, and we can do it continuously along any sweepout without changing the homotopy class, see appendix D in \cite{CM2008}. Therefore, we have
\begin{equation}\label{defWanergy} W(N,M,\omega)= \inf_{\sigma \in \omega} \max_{t\in \mathbb{B}^{k-2}} \frac{1}{2}\int_{\mathbb{D}} \vert \nabla \sigma_t \vert^2. \end{equation}
The critical maps $u : \mathbb{D} \to N$ of the energy with the constraint that $u(\partial\mathbb{D})\subset M$ are the so-called free-boundary harmonic maps, that is harmonic maps such that $\partial_{\nu} u \in (T_u M)^{\perp}$ on $\partial \mathbb{D}$. Notice that looking at the Hopf differential, a harmonic map with free boundary on the disk is automatically conformal and then minimal. This means that the set of critical points is the same considering either the energy or the area. Here is our main theorem.

\begin{theo}
\label{main0}
Let $N$ a $n$-dimensional closed regular manifold and $M$ a $m$-dimensional compact submanifold, let $\omega$ a homotopy class of sweepouts such that $W(N,M,\omega)>0 $. Then, there is a minimizing sequence of sweepouts $\sigma^n \in \omega$ such that for any sequence of parameters $t_n \in \mathbb{B}^{k-2}$ satisfying
$$ Area(u_n) \to W(N,M,\omega) \hbox{ as }  n \to +\infty $$
with $u_n = \sigma_{t_n}^n$, then, up to a subsequence, there exists $r\geq0$ (possibly branched) minimal disks with free boundary in $M$, $\theta_i : \D \rightarrow N$ and $s\geq 0$ (possibly branched) minimal spheres $\omega_j: S^2\rightarrow N$ such that
$$u_j \rightarrow \sum_{i=1}^r \theta_i +\sum_{j=1}^s \omega_j$$
in the sense of the $W^{1,2}$-bubble convergence\footnote{see section 2 for a precise definition.}. Moreover, they achieve the {\it width}, i.e.
$$W(N,M,\omega)=\sum_{i=1}^r \mathrm{Area}(\theta_i) + \sum_{j=1}^s \mathrm{Area}(\omega_j).$$
\end{theo}

{\bf Remark:}  Bubble convergence implies varifold convergence, as proved in A.3 of \cite{CM2008}.\\

This gives an existence theorem as soon as $\pi_k(N,M)\neq 0$, either there is a minimal disk with free boundary, or a minimal sphere. This existence part was already obtained by Fraser \cite{F00}. Notice that our limiting surfaces are not a priori embedded and can even possess some isolated branched points. Thanks to a min-max method using geometric measure theory tools, Li \cite{Li} and Li-Zhou \cite{LZ} proved the existence of properly embedded minimal disks in the case when $M$ is the boundary of a domain in $\mathbb{R}^n$, without convexity assumption on $M$. However, their disk does not a priori achieve the {\it width}.\\

Our conclusion is optimal in the general case, since examples where the limiting surface should be a union of disconnected spheres or disks have to occur. We cannot either expect $\mathcal{C}^0$-bubble convergence, and of course, the limit is not a priori in the same homotopy class, even in the minimization case. For instance consider a manifold which contains a minimal sphere which enclosing a singularity and which is asymptotically flat, as the Schwarzschild space. If we consider $M$, not simply connected, being far from the minimal horizon and we try to minimize the area of disk in the homotopy class that encloses the singularity and bounds a non-homotopically trivial curve then we are going to blow at least one minimal horizon and possibly a minimal free boundary disk that does not enclose the horizon.

It is a very interesting question to know if the {\it width} can be achieved by a single surface. We can easily exclude interior blow-up points assuming there is no minimal sphere in $N$ which is the case as soon as the curvature of $N$ is non-positive for instance, see corollary 8.6.3 of \cite{J}. In the general case we expect some bubbling phenomena to occur at the boundary. Nevertheless we expect the situation to be much more rigid when $M$ is the boundary of a two-convex domain since in this case the minimal disk must stay in the interior.

\begin{defi}[Definition 1.2 of \cite{F02}] A hypersurface $M$ in a Riemannian manifold $N$ is two-convex if the sum of any pair of principal curvatures of $M$ with respect to the inward pointing unit normal is positive. A two-convex domain is a domain with smooth two-convex boundary.
\end{defi}

Since, the boundary of a two-convex domain is a barrier for minimal disk, see proof of theorem 2.1 in \cite{F02}, we have the following immediate corollary

\begin{cor}
\label{main1}
Let $N$ a  two-convex domain in $\R^p$ or in manifold with non-positive curvature, let  $\omega$  a homotopy class of sweepouts such that $W(N,\partial N,\omega)>0 $
then there is a minimizing sequence of sweepouts $\sigma^n \in \omega$ such that for any sequence of parameters $t_n \in \mathbb{B}^{k-2}$ satisfying
$$ Area(u_n) \to W(N,\partial N,\omega) \hbox{ as }  n \to +\infty $$
with $u_n = \sigma_{t_n}^n$, then, up to a subsequence, there exist $r\geq1$ minimal disks with free boundary in $\partial N$, $\theta_i : \D \rightarrow N$  such that
$$u_n \rightarrow \sum_{i=1}^r \theta_i$$
in the sense of the $W^{1,2}$-bubble convergence. Moreover, they achieve the {\it width}, i.e.
$$W(N,\partial N,\omega)=\sum_{i=1}^r \mathrm{Area}(\theta_i) .$$
\end{cor}

Finally, inspired by the result of Fraser, theorem 1 of \cite{F00}, see also \cite{Ri18}, we should be able to prove in an incoming work that the sum of the indices of the minimal disks is at most $k-2$, hence using the fact that the index of a minimal disk in a convex domain of $\R^p$ is at least $[(p-2)/2]$, see theorem 2.5 of \cite{F00}, we will have an upper bound on $r$ in the convex case, especially that $r=1$ for a convex domain in $\R^3$ and get a universal bound on $r$ with respect to the topology in $\R^p$ with $p\geq 3$. For a general domain, when $r>1$ it will be a very interesting question to localize the blow-up point and to describe the possible configuration of free boundary minimal disks that can occur .\\

Our study can be viewed as the free-boundary counterpart of the classical existence results of harmonic maps on a 2-sphere into closed manifolds initiated by Sacks-Uhlenbeck \cite{SU} in the minimization case and Micallef-Moore \cite{MM} for the min-max. Even in the closed case, the main problem to apply the classical methods for min-max is that we need some Palais-Smale assumption for the energy functional. We need that a sequence of maps $u_n \in W^{1,2}(\mathbb{S}^2, N)$ satisfying
\begin{equation} \label{palaissmale} \Delta u_n = A(u_n)(\nabla u_n,\nabla u_n) + f_n \end{equation}
converges in $W^{1,2}$ up to a subsequence, where $A$ is the second fundamental form given by the embedding $N\subset \mathbb{R}^p$ and $f_n \to 0$ as $n \to +\infty$ in $W^{-1,2}$. Of course, when $f_n =0$, it is the harmonic map equation. It is well known that even for $f_n = 0$, we cannot expect $W^{1,2}$-convergence. The optimal result was given by Parker, see \cite{PAR1996}, with a $W^{1,2}\cap \mathcal{C}^0$-bubble convergence. However, in this paper, Parker also gave an example of $u_n$ satisfying \eqref{palaissmale} which cannot even converge in the $W^{1,2}$-bubble convergence sense. 

Notice that we have more information than \eqref{palaissmale} since the problem comes from a min-max, and even if we renounce to prove that any minimizing sequence for the min-max converges we can select a suitable one by a regularization process. There are two classical strategies. 

The first one, the so-called viscosity method is to change the functional depending on a smoothing parameter $\alpha$ so that it satisfies the Palais-Smale assumption for the regularized functional and prove the convergence of the critical min-max solutions associated to $\alpha$ as $\alpha \to \alpha_0$. In this case, we have of course a sequence of solutions of \eqref{palaissmale}, but they are also critical points of the regularized functional. Therefore we add another structure on the sequence of equations and we can expect convergence. This approach was introduced by Sacks-Uhlenbeck using the so-called $\alpha$-energy
$$ E_{\alpha}(u) = \frac{1}{2} \left(  \int_{\mathbb{S}^2} \left(1+\vert\nabla u\vert^2\right)^{\alpha}  - \int_{\mathbb{S}^2} 1 \right) $$
where $E_{\alpha} \to E$ as $\alpha\to 1$. It suffices to study a sequence of critical points of $E_{\alpha}$, satisfying the "$\alpha$-harmonic map equation" and hope for some convergence as $\alpha \to 1$. The main step they proved is $\eps$-regularity independent from $\alpha$, leading to the celebrated existence theorem. Fraser used the same method to prove the existence part in the free-boundary case, adding an important work for free-boundary regularity. However, we need another step in order to prove bubble convergence and energy identities : the "no-neck energy" lemma. Unfortunately, such a lemma is not true a priori for a sequence of $\alpha$-harmonic maps, this was proved by Li-Wang \cite{LW}. Again, of course, we have more information than the $\alpha$-harmonic map equation : the sequence comes from min-max solutions, and thanks to the monotonicity trick by Struwe \cite{Struwe}, up to a subsequence, one can add an entropic condition that Lamm used \cite{Lamm10} to prove the $W^{1,2}$-bubble convergence. The viscosity method was also performed by Rivi\`ere in the extrinsic case, adding to the area functional \cite{Riv17} or Willmore functional \cite{Riv15} adapted smoothers depending on a $L^p$-norm of the second fundamental form.

The second strategy is to start from a minimizing sequence and replace it with a suitable procedure by a new competitor which satisfies conditions leading to compactness. This idea goes back to Birkhoff in the case of closed geodesics, with the famous curve shortening process : choosing sets of points sufficiently close, we replace each portion between two successive points by the unique geodesic which joins these points. Then, we can use the compactness of the space of sweepouts of piecewise geodesic curves. We use this strategy to prove Theorem \ref{main}, based on the modern adaptation by Colding and Minicozzi \cite{CM2008}, and their proof of $W^{1,2}$-bubble convergence for the min-max problem for 2-spheres in a compact manifold diffeomorphic to $\mathbb{S}^3$. Generalizations for tori and higher genus surfaces was performed by  Zhou in \cite{Z10} and \cite{Z17} where harmonic maps are not automatically conformal. In the next section, we will give the main steps of proof of Theorem \ref{main} (implying Theorem \ref{main0}), which handles this strategy for free boundary issues.\\

In order to perform the replacement procedure, one needs to get at least the uniqueness of harmonic maps with small energy with respect to the Dirichlet boundary data. Even though it is a direct consequence of maximum principle when the target manifold is $\R^p$ it is much harder in the general case. Colding and Minicozzi proved in fact much more, since they proved the energy convexity of the energy functional using some property of holormorphic functions, see section C of \cite{CM2008}. In fact, the crucial point is the compactness compensation phenomena which appears in all conformally invariant problems. This remark permits the first author and Lin \cite{LL} to prove the convexity for some biharmonic maps, see also \cite{L17} for a review. In the present paper we prove the following free boundary version of Colding-Minicozzi energy convexity.

\begin{theo}\label{econvexb} Let $N$ be a compact submanifold of $\mathbb{R}^n$ and $M$ a closed submanifold of $N$. Then, there exists a constant $\eps_0 >0$  such that if $u,v \in W^{1,2}(\mathbb{D}_+,N)$ with $u\vert_{A} = v\vert_{A} $, for almost every $x\in I$, $v(x)\in M$ and $u(x)\in M$, $E(u) \leq \eps_0$ and $u$ is weakly harmonic meeting $M$ orthogonally along $I$, then we have the energy convexity
\begin{equation}
\frac{1}{2}\int_{\mathbb{D_+}} \left\vert \nabla (v - u)\right\vert^2 \leq \int_{\mathbb{D_+}} \left\vert \nabla v \right\vert^2 - \int_{\mathbb{D_+}} \left\vert \nabla u \right\vert^2 \hskip.1cm.
\end{equation}
Where $\mathbb{D}_+= \mathbb{D} \cap \{(x_1,x_2)\, \vert \, x_2> 0\}$ and $\partial \mathbb{D}_+ = A \cup I$ with $I = (-1,1)\times \{0\}$ and $A = \partial \mathbb{D}_+ \setminus I$.
\end{theo}

This kind of energy convexity is also of first importance to get some strong convergence of the flow associated to some conformally invariant problem such as harmonic maps or bi-harmonic maps, see \cite{Lin} and \cite{LL}. Moreover, using the flow as another smoother for sweepouts could lead to some compactness, as in Fraser-Schoen \cite{FS}. Notice that energy identities for harmonic flow have already been performed, see for instance \cite{DT}.\\

Organization of the paper: In section 2, we give a detailed statement of our main theorem and we give a detailed sketch of the proof, since the proof is somehow technical. The strategy is the one of Colding-Minicozzi \cite{CM2008}. Section 3 is devoted to the proof of energy convexity. In section 4, we perform  the replacement procedure to produce the energy decreasing map. Section 5, is the key point of the paper since we prove the no neck energy which will ensure that the width is achieved. Section 6 is dedicated to the proof of the main theorem. Finally in the appendix, we give our version of the $\eps$-regularity for free boundary harmonic maps.

\medskip

{\bf Added in proof :}{\it We have been told recently by Lin, that with his collaborator Zhou and Sun, they prove a similar result.}

\medskip

{\bf Acknowledgements} : 
{\it The authors are very grateful to the referees for their very careful reading of the paper.}

\section{ Main theorem and sketch of the proof}

In this section we give the full statement of our main theorem and we detail the strategy of the proof in what follows.

\begin{theo}
\label{main}
Let $N$ a $n$-dimensional closed regular manifold and $M$ a $m$-dimensional compact submanifold, then 
$$ \pi_k(N,M)\neq 0 \Leftrightarrow \exists \omega \hbox{ a homotopy class of sweepouts such that } W(N,M,\omega)>0 $$
and in this case, there is a minimizing sequence of sweepouts $\sigma^n \in \omega$ such that for any sequence of parameters $t_n \in \mathbb{B}^{k-2}$ satisfying
$$ Area(u_n) \to W(N,M,\omega) \hbox{ as }  n \to +\infty $$
with $u_n = \sigma_{t_n}^n$, then, up to a subsequence, there exists $r\geq 0$, $s\geq 0$ and
\begin{itemize}
\item $u_\infty: \mathbb{D}\to N$ (possibly constant or branched) minimal disks with free boundary on $M$, 
\item $\theta_1,\cdots,\theta_r : \mathbb{D}\to N$ non constant (possibly branched) minimal disks with free boundary on $M$, some centers $a^n_1,\cdots,a^n_r \in \partial\mathbb{D}$, and scales $\lambda_1^n,\cdots,\lambda_r^n \to 0$ as $n\to +\infty $,
\item $\omega_1,\cdots,\omega_s : \mathbb{S}^2 \to N$ non constant  (possibly branched) minimal spheres, some centers  $b^n_1,\cdots,b^n_s \in \mathbb{D}$ and scales $\nu_1^n,\cdots,\nu_s^n \to 0$ as $n\to +\infty $,
\end{itemize}
such that 
\begin{itemize}
\item $u_n$ converges to $u_{\infty}$ in $W^{1,2}_{loc}(\D\setminus\{b_1^\infty,\dots,b^\infty_s\})$,
\item $u_n \circ \phi_{a^n_i}( \lambda^n_{i} \phi_{1}^{-1}(.) )$ converges  to $\theta^i$ in $W_{loc}^{1,2}(\mathbb{D}\setminus \{S_i\})$ for $1\leq i \leq r$, where $S_i$ is a finite set,
\item $u_n(b_j^n+ \nu^n_j . )$ converges to $\omega^j \circ \pi^{-1}$ in $W_{loc}^{1,2}(\mathbb{R}^2\setminus\{S_j\})$ for $1\leq j \leq s$, where $S_j$ is a finite set,
\end{itemize}
where $\pi : \mathbb{S}^2 \to \mathbb{R}^2$ is the stereographic projection with respect to the north pole and $\phi_a : \mathbb{D} \to \mathbb{R}_+^2$ satisfies $\phi_a(z) = i \frac{a-z}{a+z}$. Moreover, we have the energy identity :
$$ \lim_{n\to +\infty} \int_{\mathbb{D}}\vert\nabla u_n \vert^2 = \int_{\mathbb{D}}\vert\nabla u_{\infty} \vert^2 + \sum_{i=1}^r  \int_{\mathbb{D}}\vert\nabla \theta^i \vert^2 +   \sum_{i=1}^s \int_{\mathbb{D}}\vert\nabla \omega^i \vert^2 \hskip.1cm.$$
\end{theo}

The conclusion of the theorem gives a classical so-called "$W^{1,2}$ bubble convergence", see \cite{PAR1996},\cite{Struwe} and \cite{LR}. Notice that the parameters satisfy
$$ \frac{\vert a_{\alpha}^n - a_{\beta}^n \vert}{\lambda_{\alpha}^n + \lambda_{\beta}^n} + \frac{\lambda_{\alpha}^n}{\lambda_{\beta}^n} + \frac{\lambda_{\beta}^n}{\lambda_{\alpha}^n} \to +\infty \hbox{ as } n\to +\infty $$
for $1\leq \alpha<\beta\leq r$ and
$$ \frac{\vert b_{\alpha}^n - b_{\beta}^n \vert}{\nu_{\alpha}^n + \nu_{\beta}^n} + \frac{\nu_{\alpha}^n}{\nu_{\beta}^n} + \frac{\nu_{\beta}^n}{\nu_{\alpha}^n} \to +\infty \hbox{ as } n\to +\infty $$
and
$$ \frac{1-\vert b_{\alpha}^n \vert}{\nu_{\alpha}^n} \to +\infty \hbox{ as } n\to +\infty$$
for $1\leq \alpha<\beta\leq s$.

\medskip

In the rest of the paragraph we detail the strategy of the proof. As we already mentioned, the proof follows the strategy of Colding-Minicozzi for the sphere case, we will mainly focus on the new difficulties given by the presence of a boundary. Especially we have to develop new tools since the classical argument which permits to get estimates by reducing the size of the set does not work at the boundary.\\

By Nash's theorem \cite{Nash}, we may assume in the rest of the paper that $N$ is isometrically embedded in some $\R^p$.\\

\textit{STEP 0: $W=W_E$}

Setting
$$W_E(N,M,\omega)= \inf_{\sigma \in \omega} \max_{t\in \mathbb{B}^{k-2}} \frac{1}{2}\int_{\mathbb{D}} \vert \nabla \sigma_t \vert^2, $$
we have $W_E=W$.\\

This was the seminal idea of Douglas-Rad\'o to work on the Energy rather than the $Area$. The main idea is that any $W^{1,2}$-map can be reparametrized as a quasi-conformal map for which the Area and the Energy are as close as desired. Colding and Minicozzi carefully improved this idea to the sweepout setting. The main point is that the reparametrization has to be continuous with respect to the parameters of the sweepout. Since the sphere case and the disk case are similar, they both possess only one conformal class, we just remind the main three steps of appendix D of \cite{CM2008}.\\

\begin{itemize}
\item It is clear that $W\leq W_E$. Then let us consider $\sigma$ such that $\ds \max_{t\in \mathbb{B}^{k-2}} Area(\sigma_t) \leq W +\frac{\eps}{2}$, for some $\eps>0$.
\item Then, in lemma D.1 of \cite{CM2008}, improving  the density argument of Schoen-Uhlenbeck, see last proposition of \cite{ScU}, Colding-Minicozzi explained how to regularize the sweepout $\sigma$ to some $\tilde{\sigma}\in C^0(\mathbb{B}^k,C^2(\D,N))$, verifying $\ds \max_{t\in \mathbb{B}^{k-2}} Area(\tilde{\sigma}_t) \leq W +\eps$. We have to notice that in our case the boundary of the mollified sweepout is {\it a priori} not in $M$. In order to solve this issue, it suffices to consider the family of $C^0$-curves $c_t={\sigma_t}_{\vert \partial \D}$, and to regularize it to some $\tilde{c}_t$ which still take values into $M$. Then we consider $\tilde{v}_t=\pi_N(h_t)$ where $h_t$ is the harmonic extension of $\tilde{c}_t$ into $\R^p$ and $\pi_M$ the projection onto $N$. This $\tilde{v}_t$  is a smooth map close to $\tilde{\sigma}_t$ in some neighborhood of $\partial \D$. Finally it suffices to interpolate $\tilde{\sigma}_t$ and $\tilde{v}_t$ in this neighbourhood to get a $C^2$ sweepout, still denotes $\tilde{\sigma}_t$, which is $W^{1,2}\cap C^0$ closed to $\sigma_t$ and which sends $\partial \D$ to $M$.\\
An alternative approach for the "projection" of $\tilde{\sigma}$ consists in considering, for $\delta>0$, a map $F_\delta:N\rightarrow N$ such that $\vert DF_\delta\vert =1 +O(\delta)$ such as ${F_\delta}_{\vert N \setminus M_\delta}=Id$, see section 4 for the definition of $M_\delta$, and $F_\delta$ is a retract of $M_{\delta^2}$ onto $M$.\\
\item Considering the pullback metric $\tilde{\sigma}_t^*(\xi)$ on the disk, it can be degenerated, hence one considers $g_{\delta,t}=\tilde{\sigma}_t^*(\xi) +\delta \vert dz\vert^2$ with $\delta>0$. 
Then thanks to the Riemann mapping theorem for variable metric of Ahlfors-Bers \cite{AB}, see also section 3.2 of Jost \cite{Jo}, one finds a unique conformal diffeomorphism $h_{\delta,t}$ which fixes  three points on the boundary and pulls back $g_{\delta,t}$  to $\vert dz\vert^2$, with some $W^{1,2}\cap C^0$ control on the diffeomorphism.

\item Finally using this control and choosing $\delta>0$ small enough, Colding-Minicozzi proved that the Energy of $\tilde{\sigma}_t\circ h_{\delta,t}$ can be made as close as we want to the $Area$ of $\sigma_t$.
\end{itemize}

\medskip

\textit{STEP 1: $\eps$-regularity for a sequence of harmonic maps with free boundary of uniformly bounded energy}

\medskip

An $\eps$-regularity convergence for a sequence of critical points has to be true for Theorem \ref{main}. This result was proved by Laurain-Petrides \cite{LP} when $N=\mathbb{B}^3$ and $M=\mathbb{S}^2$ thanks to reflexion methods and the use of Rivi\`ere's conformal invariant equations theory. Jost-Liu-Zhou generalized this result for any $N$ and $M$, see \cite{JLM}. For sake of completeness, we recall the $\eps$-regularity result in the appendix (see Proposition \ref{epsreg}) with a concise argument.

\medskip

\textit{STEP 2: Convexity for free boundary energy}

\medskip

This is a key step of the paper of Colding and Minicozzi in order to define a replacement procedure. In the classical problem of geodesics by Birkhoff, we replace portions of curves by a geodesic which joins their ends, but we need uniqueness of geodesics: the points have to be below the injectivity radius. We state here an analogous result for surfaces.

Let $\mathbb{D}_+= \mathbb{D} \cap \{(x_1,x_2)\, \vert \, x_2> 0\}$. We denote $\partial \mathbb{D}_+ = A \cup I$ where $I = (-1,1)\times \{0\}$ and $A = \partial \mathbb{D}_+ \setminus I$. The next result states that the energy functional is strictly convex around small energy maps with the constraint that the boundary $I$ is sent in some submanifold $M$ of $N$.

\begin{theo} Let $N$ be a compact submanifold of $\mathbb{R}^p$ and $M$ a closed submanifold of $N$. Then, there exists a constant $\eps_0 >0$  such that if $u,v \in W^{1,2}(\mathbb{D}_+,N)$ with $u\vert_{A} = v\vert_{A} $, for almost every $x\in I$, $v(x)\in M$ and $u(x)\in M$, $E(u) \leq \eps_0$ and $u$ is weakly harmonic meeting $M$ orthogonally along $I$, then we have the energy convexity
\begin{equation}
\frac{1}{2}\int_{\mathbb{D_+}} \left\vert \nabla (v - u)\right\vert^2 \leq \int_{\mathbb{D_+}} \left\vert \nabla v \right\vert^2 - \int_{\mathbb{D_+}} \left\vert \nabla u \right\vert^2 \hskip.1cm.
\end{equation}
\end{theo}

We can deduce the following uniqueness result as an immediate corollary.

\begin{cor} Let $N$ be a compact submanifold of $\mathbb{R}^p$ and $M$ a closed submanifold of $N$. Then, there exists a constant $\eps_0 >0$  such that if $u,v \in W^{1,2}(\mathbb{D}_+,N)$ are weakly harmonic meeting $M$ orthogonally along $I$, with $u\vert_{A} = v\vert_{A} $, $E(u)\leq \eps_0$ and $E(v) \leq \eps_0$ then $ u\equiv v$.
\end{cor}

Notice that the previous theorem is  a "free-boundary" version of the energy convexity result by Colding and Minicozzi in \cite{CM2008} :

\begin{theo}\label{econvex} Let $N$ be a compact submanifold of $\mathbb{R}^n$. Then, there exists a constant $\eps_0 >0$  such that if $u,v \in W^{1,2}(\mathbb{D},N)$ with $u\vert_{\partial \mathbb{D}}  = v\vert_{\partial \mathbb{D}} $, $E(u) \leq \eps_0$ and $u$ is weakly harmonic, then we have the energy convexity
\begin{equation}
\frac{1}{2}\int_{\mathbb{D}} \left\vert \nabla (v - u)\right\vert^2 \leq \int_{\mathbb{D}} \left\vert \nabla v \right\vert^2 - \int_{\mathbb{D}} \left\vert \nabla u \right\vert^2 \hskip.1cm.
\end{equation}
\end{theo}

Besides the uniqueness consequence, this theorem is very useful for many other problems such as flow convergence, see \cite{LL}. In fact our proof simplifies the original proof by Colding-Minicozzi and the one of Lamm-Lin \cite{LaL}. In fact our idea applies probably to all conformally invariant problem, since it relies only on the $\eps$-regularity, the fact that the right hand side of the equation is orthogonal to $TN$. Those two facts were the key ingredients of the proof of a quantization phenomena for all conformally invariant problem, see \cite{LR}. The last ingredient is a Hardy inequality. In order to get the convexity for the free boundary setting we need to generalize the Hardy inequality to free boundary problems. It is something already known, see theorem 1 of \cite{FMT}, but for the sake of completeness we give a complete statement and a full proof see theorem \ref{HB}. We would like to remark that the Hardy inequality can be seen as a special case of some weighted Poincar\'e inequality, or even Poincar\'e inequality on non-compact manifolds. We will discuss this point of view in an incoming paper \cite{LLP}.

\medskip

\textit{STEP 3: The harmonic replacement procedure and decreasing energy map}

\medskip

In the proof of Theorem \ref{main}, we first take a minimizing sequence of sweepouts. Of course, we extract a classical Palais-Smale sequence $u_n$ (satisfying \eqref{palaissmale} in the interior) but as already said, we cannot conclude for $W^{1,2}$-bubble convergence. We aim at replacing this sequence by another one which will satisfy stronger Palais-Smale-like properties used in \textit{STEP 4} to prove convergence. This is stated in Theorem \ref{harmrepl}. Let's sketch this procedure in the particular case of minimization (the set of parameters is trivial $\mathbb{B}^{k-2} = \{0\}$).

Let $u \in \mathcal{A}$. Let $\mathcal{E}$ the set of finite families of disjoint element of $\mathfrak{B}$ and $\mathfrak{HB}$ where 
$$\mathfrak{B}=\{ \overline{B(a,r)}\hbox{ s.t. } B(a,r)\subset \mathbb{D}\}$$
and 
$$\mathfrak{HB}=\{\overline{B(a,r)\cap \mathbb{D}} \hbox{ s.t. } a\in \partial\mathbb{D} \hbox{ and } \partial B(a,r) \hbox{ intersect orthogonally}\}.$$

For $\alpha\in (0,1]$, and $\mathcal{B}\in \mathcal{E}$ we denote by $\alpha\mathcal{B}$ the collection of disks concentric to those of $\mathcal{B}$ with radius dilated by $\alpha$ and half disks which meet the boundary orthogonally, with center collinear to  those of $\mathcal{B}$  and with radius dilated by $\alpha$.

For and $\mathcal{B} \in \mathcal{E}$ such that the energy of $u$ on $\bigcup_{B\in \mathcal{B}} B$ is less than the  $\eps_0$ of Theorem \ref{econvexb}, then we denote by $H(u,\mathcal{B}):\mathbb{D}\rightarrow N$ the map that coincides 
\begin{itemize}
\item on $\mathbb{D}\setminus \bigcup_{B\in \mathcal{B}} B$ with $u$
\item on $\bigcup_{B\in \mathcal{B}} B$ with the (unique) energy minimizing map from $\bigcup_{B\in \mathcal{B}} B$ to $N$ that agrees with $u$ on $\bigcup_{B\in \mathcal{B}} \partial B \setminus \partial{\mathbb{D}}$ with the constraint that it lies in $M$ on $\bigcup_{B\in \mathcal{B}} \partial B \cap \partial{\mathbb{D}}$.
\end{itemize}

We aim at decreasing the energy of $u$ as much as we can. We set
$$ e_u = \sup \left\{ E(u) - E(H(u,\mathcal{B})) ; \mathcal{B}\in\mathcal{E} , \int_{\bigcup_{B\in \mathcal{B}} B} \vert \nabla u\vert^2 \leq \eps_0 \right\} \hskip.1cm,$$
so that if $u$ is not already harmonic, we can pick some $\widetilde{\mathcal{B}}$ such that 
\begin{equation} \label{bigdecrenergy} E(u) - E(\tilde{u}) \geq \frac{e_u}{2}  \end{equation}
where $\tilde{u} = H(u,\widetilde{\mathcal{B}})$ is the replaced map we choose. For this map, we can easily prove thanks to Theorem \ref{econvexb} and Theorem \ref{econvex} in \textit{STEP 2}, an exchange formula (see \eqref{exchange2} in lemma \ref{exchangeharm}) and the definition of $e_u$ that for any $\mathcal{B}$

\begin{eqnarray*}
\frac{1}{4}\int_{\mathbb{D}} \vert \nabla \tilde{u} - \nabla H(\tilde{u},\frac{1}{2}\mathcal{B})\vert^2 & \leq & E(\tilde{u}) - E(H(\tilde{u},\frac{1}{2}\mathcal{B})) \\
& \leq &  E(u) - E(H(u,\mathcal{B})) + \frac{\left( E(u)- E(\tilde{u}) \right)^{\frac{1}{2}}}{\kappa} \\
& \leq & e_u + \frac{\left( E(u)- E(\tilde{u}) \right)^{\frac{1}{2}}}{\kappa} 
\end{eqnarray*}
and with the definition of $\tilde{u}$ (see \eqref{bigdecrenergy}) we conclude that
\begin{equation} \label{epsregesttildeu} \int_{\mathbb{D}} \vert \nabla \tilde{u} - \nabla H(\tilde{u},\frac{1}{2}\mathcal{B})\vert^2 \leq C \left( E(u)- E(\tilde{u}) \right)^{\frac{1}{2}} \hskip.1cm.\end{equation}
for any $\mathcal{B} \in \mathcal{E}$ with $\int_{\bigcup_{B\in \mathcal{B}} B} \vert \nabla u\vert^2 \leq \eps_0$, up to decrease $\eps_0$ and for some constant $C$ which only depends on $M$ and $N$.

Therefore, if $\{u_n\}$ is a minimizing sequence (we recall that we assumed that the set of parameters is trivial $\mathbb{B}^{k-2} = \{0\}$), we can define a new sequence as previously $\{\tilde{u}_n\}$ satisfying \ref{epsregesttildeu}. By construction, $\{u_n\} \in \mathcal{A}$ and $E(\tilde{u}_n) \leq E(u_n)$ and by continuity of the harmonic replacement (see proposition \ref{continuity}), by shrinking the radii of the balls and half balls used to define $\tilde{u}_n$, we see that $u_n$ and $\tilde{u}_n$ lie in the same homotopy class. Therefore, $\{\tilde{u}_n\}$ is also a minimizing sequence so that thanks to \ref{epsregesttildeu},
$$ \int_{\mathbb{D}} \vert \nabla \tilde{u}_n - \nabla H(\tilde{u}_n,\frac{1}{2}\mathcal{B})\vert^2 \to 0 \hbox{ as } n\to +\infty $$
for any $\mathcal{B} \in \mathcal{E}$ such that $\int_{\bigcup_{B\in \mathcal{B}} B} \vert \nabla u_n\vert^2 \leq \eps_0$. 

As already said, the general case is given in Theorem \ref{harmrepl}. The idea is the same but much more technical since we have to define a replacement procedure all along the sweepouts. This is possible by an intensive use of proposition \ref{continuity}.

\medskip

\textit{STEP 4: A specific Palais-Smale assumption}

\medskip 

Thanks to the previous steps, we can build a minimizing sequence of sweepouts $\sigma^n \in \omega$ such that for any sequence of parameters $t_n \in \mathbb{B}^{k-2}$ satisfying
$$ Area( u_n ) \to W(N,M,\omega) \hbox{ as }  n \to +\infty $$
with $u_n = \sigma_{t_n}^n$, we have the assumptions \eqref{epsregesttildeuj} and \eqref{quasiconfuj} of the following theorem implicitly proved in section 5

\begin{theo} \label{palaissmaleassump} Let $u_n\in \mathcal{A}$ be a sequence of maps with uniformly bounded energy such that
\begin{equation} \label{epsregesttildeuj} \int_{\mathbb{D}} \vert \nabla u_n - \nabla H(u_n,\eta\mathcal{B})\vert^2 \to 0 \hbox{ as } n\to +\infty \end{equation}
for any $\mathcal{B} \in \mathcal{E}$ such that $\int_{\bigcup_{B\in \mathcal{B}} B} \vert \nabla u_n\vert^2 \leq \eps_0$ and
\begin{equation} \label{quasiconfuj} Area(u_n) = E(u_n) + o(1) \hbox{ as } n\to +\infty \hskip.1cm. \end{equation}
Then, up to a subsequence, $\{u_n\}$ $W^{1,2}$-bubble converges.
\end{theo}

We can consider \eqref{epsregesttildeuj} as a kind of $\eps$-regularity property for minimizing sequences of the energy. Indeed, outside a finite number of points where the energy could concentrate over $\eps_0$, we deduce a strong convergence in $W^{1,2}$ to a harmonic map. \eqref{quasiconfuj} is a crucial assumption in order to get a no-neck energy lemma (see proposition \ref{noneckenergy}). Usually, for harmonic equations, we use a Poho\v{z}aev identity to prove a no-neck energy lemma. Here, we do not have any equation but \eqref{quasiconfuj} is a quasi-conformal assumption, where being conformal is even stronger than a Poho\v{z}aev identity. It is why we can generalize Theorem \ref{palaissmaleassump} even in the case where $u_n$ is defined on a set of degenerating conformal classes on a surface (this work was performed by Zhou in the closed case).

\medskip

\section{ Proof of theorem \ref{econvexb} }

Before proving theorem \ref{econvexb}, we have to prove some generalization of the Hardy inequality. Before that, let us remind the classical case, see section 1.3.1 of \cite{Maz}, see also \cite{MMP}.

\begin{theo}
\label{HA} Let $u\in W^{1,2}_0(\mathbb{D})$ then
\begin{equation} 
\frac{1}{4} \int_D \frac{u^2}{(1-\vert x\vert)^2}dx \leq  \int_{\mathbb{D}} \left\vert\nabla u \right\vert^2 dx \hskip.1cm.
\end{equation}
\end{theo}

For our purpose we prove an trace version of this classical result, largely inspired by theorem 1.1 in \cite{FMT}.

\begin{theo}
\label{HB} Let $u\in W^{1,2}(\mathbb{D}_+)$ such that the trace of $u=0$ on $A$. Then, the trace of $u$ on $I$ which we still denote $u$ belongs to $L^2\left(I, \frac{1}{1-x_1^2}\right)$ and
\begin{equation} \label{THS}
\int_I \frac{u^2}{1-x_1^2}dx_1 \leq \frac{\pi}{2} \int_{\mathbb{D}_+} \left\vert\nabla u \right\vert^2 dx_1dx_2 \hskip.1cm.
\end{equation}

\end{theo}

We would like to remark that the Hardy inequality can be seen as a special case of some Weighted Poincar\'e inequality, or even Poincar\'e inequality on a non-compact manifold. We will discuss this point of view in an incoming paper \cite{LLP}.\\

{\it Proof of theorem \ref{HB} :}\\

We let $\phi\in \mathcal{C}^{\infty}(I\times \mathbb{R})$ a positive function and $u \in \mathcal{C}_c^{\infty}(I\times \mathbb{R})$. Then
$$ \int_{I\times [0,+\infty]} \left\vert \nabla u - \frac{\nabla \phi}{\phi} u \right\vert^2 = \int_{I\times [0,+\infty]} \left\vert \nabla u \right\vert^2 + \frac{\left\vert \nabla \phi \right\vert^2 }{\phi^2}u^2 - \int_{I\times [0,+\infty]} \frac{\nabla\phi}{\phi} \nabla (u^2) \hskip.1cm.$$
After an integration by parts, we get
$$ - \int_I \frac{\phi_{x_2}}{\phi} u^2= \int_{I\times [0,+\infty]} \left\vert \nabla u \right\vert^2 + \int_{I\times [0,+\infty]} \frac{\Delta \phi}{\phi} u^2 - \int_{I\times [0,+\infty]} \left\vert \nabla u - \frac{\nabla \phi}{\phi} u \right\vert^2 \hskip.1cm.$$
Setting
$$ \phi(x_1,x_2) = 1 - \frac{2}{\pi} \arctan\left( \frac{x_2}{1-x_1} \right) \hskip.1cm,$$
we get that 
$$ - \frac{\phi_{x_2}}{\phi}(x_1,0) = \frac{2}{\pi}\frac{1}{1-x_1} $$
and that $\Delta \phi =0$ so that
\begin{equation} \label{THS1} \int_I  \frac{u^2}{1-x_1} \leq \frac{\pi}{2} \int_{I\times [0,+\infty]} \left\vert \nabla u \right\vert^2  \hskip.1cm. \end{equation}
By an analogous computation, letting $\phi(x_1,x_2) = 1 - \frac{2}{\pi} \arctan\left( \frac{x_2}{1+x_1} \right)$ we get
\begin{equation} \label{THS2} \int_I  \frac{u^2}{1+x_1} \leq \frac{\pi}{2} \int_{I\times [0,+\infty]} \left\vert \nabla u \right\vert^2  \hskip.1cm. \end{equation}
Summing \eqref{THS1} and \eqref{THS2} gives  \eqref{THS} for $u \in \mathcal{C}_c^{\infty}(\mathbb{D}_+\cup I)$.\\

To conclude, let $u\in W^{1,2}(\mathbb{D}_+)$ such that the trace of $u=0$ on $A$. By density, let $u_n \in \mathcal{C}_c^{\infty}(\mathbb{D}_+\cup I)$ be a sequence such that $u_n$ converges to $u$ in $W^{1,2}(\mathbb{D}_+)$. In particular up to a subsequence, $u_n$ converges to the trace of $u$ on $I$ almost everywhere. By Fatou theorem, 
$$ \int_I \frac{u^2}{1-x_1^2} \leq \liminf_{n\to +\infty} \int_I \frac{u_n^2}{1-x_1^2} \leq \liminf_{n\to +\infty} \frac{\pi}{2} \int_{\mathbb{D}_+} \left\vert\nabla u_n \right\vert^2  = \frac{\pi}{2} \int_{\mathbb{D}_+} \left\vert\nabla u \right\vert^2  \hskip.1cm, $$
which completes the proof of Theorem \ref{HB}.

\hfill $\diamondsuit$

\medskip

{\it Proof of theorem \ref{econvexb}:}\\

Let $u,v \in W^{1,2}(\mathbb{D}_+,N)$ with $u\vert_{A} = v\vert_{A} $, for almost every $x\in I$, $v(x)\in M$ and $u(x)\in M$, $E(u) \leq \eps_0$ and $u$ is harmonic meeting $M$ orthogonally along $I$. Notice that by proposition \ref{epsreg}, $u$ is a smooth map until the boundary $I$. Then, setting

\begin{equation} \label{deftheta}
\Theta = \int_{\mathbb{D_+}} \left\vert \nabla v \right\vert^2 - \int_{\mathbb{D_+}} \left\vert \nabla u \right\vert^2 - \int_{\mathbb{D_+}} \left\vert \nabla (v - u)\right\vert^2
\end{equation}
we get
\begin{equation} \label{thetaintbyparts}
\Theta  =   2 \int_{\mathbb{D_+}} \left\langle \nabla (v - u),\nabla u \right\rangle 
 =  2 \int_{\mathbb{D_+}} \left( v -  u \right) . \Delta u + 2 \int_{I}  \left(v - u\right).\partial_{\nu} u 
\end{equation}
by an integration by parts since $u=v$ on $A$. Now, we use proposition \ref{ColMinortho}, see below, twice. First, since $\Delta u = A(u)\left(\nabla u, \nabla u \right) \perp T_u N $, where $A$ is the second fundamental form of the submanifold $N$ of $R^n$, we get
\begin{equation} \label{ortho1}
\int_{\mathbb{D_+}} \left( v -  u \right) . \Delta u \geq - \left\|A\right\|_{\infty} \int_{\mathbb{D_+}} \left|\left( v -  u \right) ^{\perp} \right| \left\vert \nabla u \right\vert^2 \geq - \left\|A\right\|_{\infty} \kappa_N \int_{\mathbb{D_+}} \left| v -  u  \right|^2  \left\vert \nabla u \right\vert^2 \hskip.1cm.
\end{equation}
Also, since  $ \partial_{\nu} u \perp T_u M$ on $M$, submanifold of $\mathbb{R}^n$, we get
\begin{equation} \label{ortho2}
\int_{I} \left( v -  u \right) . \partial_{\nu} u \geq -  \int_{I} \left|\left( v -  u \right) ^{\perp} \right| \left\vert \partial_{\nu} u \right\vert \geq -  \kappa_M \int_{I} \left| v -  u  \right|^2  \left\vert \partial_{\nu} u \right\vert \hskip.1cm.
\end{equation}
Gathering \eqref{deftheta}, \eqref{thetaintbyparts}, \eqref{ortho1} and \eqref{ortho2}, we get a constant $C>0$ such that
\begin{equation} \label{thetaineq}
\Theta  \geq  - C \left( \int_{\mathbb{D_+}} \left| v -  u  \right|^2  \left\vert \nabla u \right\vert^2 + \int_{I} \left| v -  u  \right|^2  \left\vert \partial_{\nu} u \right\vert \right)
\end{equation}

which conclude the proof of theorem \ref{econvexb}, applying $\eps$-regularity which gives 
$$ \vert \nabla u \vert \leq \frac{C}{1-\vert x\vert} \sqrt{\eps_0}, $$
see theorem \ref{epsreg} in the appendix, theorem \ref{HA} and theorem \ref{HB}.

\hfill$\diamondsuit$

\medskip

We recall now a lemma by Colding and Minicozzi used to prove theorem \ref{econvexb} (see (C.11) of \cite{CM2008}).

\begin{prop}\label{ColMinortho} Let $M$ be a submanifold of $\mathbb{R}^n$. Then, there is a constant $\kappa_M >0$ such that
$$\forall p,q \in M,  \left|\left( p -  q \right) ^{\perp} \right|  \leq \kappa_M   \left| p -  q  \right|^2 $$
where $\left( p -  q \right) ^{\perp}$ denotes the normal component of $p-q$ with respect to $p\in M$.
\end{prop}

We leave to the reader the proof or Theorem \ref{econvex} thanks to proposition \ref{ColMinortho} together with the $\eps$-regularity for harmonic maps and the Hardy inequality (Theorem \ref{HA}).

\section{ Replacement procedure and Energy decreasing map }

In this section, we aim at building a systematic replacement procedure in order to regularize a minimizing sequence of sweepouts for the min-max problem, in order to obtain a sequence which satisfies a kind of Palais-Smale assumption. We precise it in Theorem \ref{harmrepl}.

In the following, we let $\delta>0$ be such that for the open neighborhoods
$$ N_{\delta} = \{ x\in\mathbb{R}^p ; d(x,N) <\delta \} \hbox{ and } M_{\delta} = \{ x\in\mathbb{R}^p ; d(x,M) <\delta \}$$ 
of the submanifolds $N$ and $M$ of $\mathbb{R}^p$ respectively, there are smooth projection maps $\pi_N : N_{\delta} \to N$ on $N$ and  $\pi_M : M_{\delta} \to M$ on $M$ such that

\begin{equation}
\label{dpn}
 \sup_{N_\delta} \Vert D \pi_N \Vert \leq 2 \hbox{ and } \forall x \in N_{\delta}, \Vert D \pi_N (x) \Vert \leq 1 + C\vert x- \pi_N(x)\vert
\end{equation}
and
\begin{equation}
\label{dpm}
 \sup_{M_\delta} \Vert D \pi_M \Vert \leq 2 \hskip.1cm.
\end{equation}

Moreover, we can extend $\pi_M : M_{\frac{\delta}{2}} \to M$ by a smooth map $\widetilde{\pi_M} : \mathbb{R}^p \to \mathbb{R}^p$. For instance, thanks to a smooth cut-off function $\chi$ such that $\chi = 1$ on $M_{\frac{\delta}{2}}$ and $\chi = 0$ on $\mathbb{R}^p \setminus M_{\delta}$, we can set $\widetilde{\pi_M} = \chi \pi_M$.

\subsection{Continuity of the harmonic extension map}

The main goal of this paragraph is to prove that the harmonic replacement is a continuous map in the free boundary case, see section 3.2 of \cite{CM2008} for the interior case. The free boundary case is far from being a simple adaptation of Colding and Minicozzi argument since we have to pay a particular attention to the boundary, especially when we construct some competitors.

\begin{prop} \label{continuity} There exists $\eps_0 >0$, such that for  every $u \in W^{1,2}\cap\mathcal{C}^0(\overline{\mathbb{D}_+},N)$ such that $u(I)\subset M$ with energy less than $\eps_0$ there is a unique energy minimizing map from $\overline{\mathbb{D}_+}$ to $N$ that agrees with $u$ on $A$ with the constraint that it lies in $M$ on $I$ denoted by $\tilde{u}$. Moreover, 
$ u \mapsto \tilde{u} $ is continuous for the $W^{1,2}\cap\mathcal{C}^0$ topology.
\end{prop}

{\it Proof of proposition \ref{continuity} :}\\

Existence is standard (see theorem 2 of 4.6 in \cite{DHS}) and uniqueness is given by Theorem \ref{econvexb}. For the continuity of the map $u\mapsto \tilde{u}$, we follow three steps.

\vspace{4mm}

\textit{STEP 1 :} There exists $C>0$ such that for any $u$ and $v$ in $W^{1,2}\cap \mathcal{C}^0(\overline{\mathbb{D}},N)$ with $u(I)\subset M$ and $v(I)\subset M$ of energy less than $\eps_0$, we have
\begin{equation} \label{lipen} \vert E(\tilde{u}) - E(\tilde{v}) \vert \leq C\left(\Vert u-v \Vert_{\mathcal{C}^0} + \Vert \nabla(u-v)\Vert_{L^2}\right). \end{equation}

First of all, up to take $C$ big enough, we can assume that $\Vert u-v \Vert_{\mathcal{C}^0} < \frac{\delta}{2 + \sup \Vert id - D\widetilde{\pi_M} \Vert}$, since $\vert E(\tilde{u}) - E(\tilde{v}) \vert \leq \eps_0$.\\

Then, we set
$$\bar{v} = \tilde{u} + \widetilde{\pi_M}(\tilde{u}+v-u) - \widetilde{\pi_M}(\tilde{u}) + v - \widetilde{\pi_M}(v) - (u-\widetilde{\pi_M}(u)) $$
in order to get some comparison function close to $\tilde{u}$ as $u$ is close to $v$ such that 
\begin{equation} \label{admass} \forall x \in I, \bar{v}(x)  \in M \hbox{ and } \forall x \in A, \bar{v}(x) = v(x) \hskip.1cm. \end{equation}
Indeed, since $\Vert u-v \Vert_{\mathcal{C}^0}\leq \frac{\delta}{2 }$, $\bar{v} = \widetilde{\pi_M}(\tilde{u}+v-u) = \pi_M(\tilde{u}+v-u) \in M $ on $I$.\\

Now, we can set
$$ \hat{v} = \pi_N(v) $$
since 
$$ d(\bar{v},N) \leq \vert \bar{v} - \tilde{u} \vert \leq \vert \widetilde{\pi_M}(\tilde{u}+v-u) - \widetilde{\pi_M}(\tilde{u}) \vert + \vert (id-\widetilde{\pi_M})(v) - (id-\widetilde{\pi_M})(u) \vert$$
so that
\begin{equation} \label{distvbar} d(\bar{v},N) \leq \left(\sup\Vert D\pi_M \Vert + \sup \Vert D(id-\widetilde{\pi_M})\Vert \right) \Vert u-v \Vert_{\mathcal{C}^0} \hskip.1cm. \end{equation}
and $d(\bar{v},N) < \delta $ by \eqref{dpm}. Then $\hat{v}$ is an admissible function for the variational characterization of $\tilde{v}$ since $\hat{v} \in N$ and $\hat{v}$ satisfies the same properties of \eqref{admass} as $\bar{v}$, hence 
\begin{equation} \label{tildevhatv} 
E(\tilde{v}) \leq E(\hat{v}). 
\end{equation}
Now, thanks to \eqref{dpn} and \eqref{distvbar},
\begin{equation} \label{hatvbarv} \int_{\mathbb{D}_+} \vert \nabla \hat{v}\vert^2 \leq \left(1 + C \Vert v-u \Vert_{C^0(\mathbb{D}_+)} \right)^2 \int_{\mathbb{D}_+} \vert \nabla \bar{v}\vert^2\hskip.1cm. \end{equation}
Now, we estimate the energy of $\bar{v}$. We have
\begin{eqnarray*}
\nabla \bar{v} - \nabla \tilde{u} & = & \left(D\widetilde{\pi_M}(\tilde{u}+v-u) - D\widetilde{\pi_M}(\tilde{u})\right).\nabla\tilde{u} \\
& &+ D\widetilde{\pi_M}(\tilde{u}+v-u).\nabla(v-u) \\
& &+ \left((id-D\widetilde{\pi_M})(v)- (id-D\widetilde{\pi_M})(u)\right) . \nabla v \\
& &+ (id-D\widetilde{\pi_M})(u). \nabla(v-u) \hskip.1cm,
\end{eqnarray*} 
so that
\begin{eqnarray*}
 \vert \nabla \bar{v} - \nabla\tilde{u} \vert & \leq & \sup \Vert D^2\widetilde{\pi_M}\Vert (\vert \nabla\tilde{u} \vert + \vert \nabla v \vert)\Vert v-u \Vert_{C^0(\mathbb{D}_+)} \\
 & & + \left( \sup \Vert D\widetilde{\pi_M} \Vert + \sup \Vert id - D\widetilde{\pi_M} \Vert \right) \vert \nabla (v - u) \vert \hskip.1cm,
\end{eqnarray*}
and there is a constant $C$ such that
\begin{equation} \label{barvtildeu}  
\int_{\mathbb{D}_+} \vert \nabla \bar{v}\vert^2 \leq \int_{\mathbb{D}_+} \vert \nabla \tilde{u}\vert^2 + C \left( \Vert v-u \Vert_{C^0(\mathbb{D}_+)} + \left( \int_{\mathbb{D}_+} \vert \nabla(v-u) \vert^2 \right)^{\frac{1}{2}} \right) \hskip.1cm.
\end{equation}
Finally, by symmetry we can assume that $E(\tilde{u}) \leq E(\tilde{v})$ so that by \eqref{tildevhatv}, \eqref{hatvbarv} and \eqref{barvtildeu}, we get STEP 1.

\vspace{4mm}

\textit{STEP 2 :} We let $u_n\to u$ in $W^{1,2}\cap\mathcal{C}^0(\overline{\mathbb{D}_+},N)$. We aim at proving that $\tilde{u}_n\to \tilde{u}$ in $W^{1,2}(\mathbb{D}_+)$.

As in STEP 1, we can set
$$ w_n = \pi_N\left( \tilde{u} + \widetilde{\pi_M}(\tilde{u}+u_n-u) - \widetilde{\pi_M}(\tilde{u}) + u_n - \widetilde{\pi_M}(u_n) - (u-\widetilde{\pi_M}(u)) \right) \hskip.1cm.$$
Adapting the proof of STEP 1, we easily prove that 
\begin{equation} \label{wntotildeu} \Vert w_n - \tilde{u} \Vert_{W^{1,2}} = o(1) \hbox{ as } n\to +\infty \hskip.1cm. \end{equation}
By STEP 1, we know that $E(\tilde{u}_n )- E(\tilde{u}) = o(1)$ as $n\to +\infty$, so that by \eqref{wntotildeu}, $E(\tilde{u}_n )- E(w_n) = o(1)$ as $n\to +\infty$. By the energy convexity, Theorem \ref{econvexb}, we deduce that $ \Vert \nabla (w_n - \tilde{u}_n) \Vert_{L^2(\mathbb{D}_+)} = o(1) $ as $n\to +\infty$. Thanks to the Poincar\'e inequality, we also have
\begin{equation}\label{wntildeun} \Vert w_n - \tilde{u}_n \Vert_{W^{1,2}(\mathbb{D}_+)} = o(1) \hbox{ as } n\to +\infty  \hskip.1cm. \end{equation}
Finally, we deduce STEP 2 from \eqref{wntildeun} and \eqref{wntotildeu}.

\vspace{4mm}

\textit{STEP 3 :} We let $u_n\to u$ in $W^{1,2}\cap\mathcal{C}^0(\overline{\mathbb{D}_+},N)$. We aim at proving that $\tilde{u}_n \to \tilde{u}$ in $\mathcal{C}^0(\overline{\mathbb{D}_+})$.

The proof of \textit{STEP 3} is written in the spirit of Qing \cite{Q}.  We crucially use that there is no concentration of energy on the boundary $A$ (which is a consequence of STEP 2) and $\eps$-regularity on symmetrized function $\tilde{u}_n$ on $\mathbb{D}$ given in the proof of Proposition \ref{epsreg}. Notice that in our case, we have to deal with the boundary $I$. We aim at proving uniform equicontinuity on $\{\tilde{u}_n\}$.

By contradiction, let $\eta >0$ sequences $x_n^1,x_n^2 \in \mathbb{D}_+$ such that up to some subsequence, 
$$ \left\vert x_n^1 - x_n^2 \right\vert \leq \frac{1}{n} \hbox{ and } \left\vert \tilde{u}_n(x_n^1) - \tilde{u}_n(x_n^2) \right\vert \geq \eta $$
and $\{x_n\}$ and $\{y_n\}$ converge to some point $x\in\overline{\mathbb{D}_+}$. We set 
$$ \delta_n^i = 1-\left\vert x_n^i\right\vert \hskip.1cm.$$
and we assume that up to a subsequence, $\delta_n^i >0$ for any $n$ (the case $\delta_n^i = 0$ is easier since the construction below is made to link $x^i_n$ to the boundary).

Thanks to Proposition \ref{epsreg}, if $1-\vert x\vert>0$, we get
$$ \eta <  \left\vert \tilde{u}_n(x_n^1) - \tilde{u}_n(x_n^2) \right\vert \leq \| \nabla \tilde{u}_n\|_{L^{\infty}\left(\mathbb{D}_{1-\frac{1-\vert x\vert}{2}}\right)} \left\vert x_n^1-x_n^2\right\vert \leq \frac{2C\left(\int_{\mathbb{D}_+} \left\vert \nabla\tilde{u}_n\right\vert^2\right)^{\frac{1}{2}}}{n (1-\vert x\vert)} \hskip.1cm.$$
Then we get a contradiction, hence the distance $\delta_n^i \rightarrow 0$ and $\int_{\mathbb{D}_{\alpha \delta_n^i}(x_n^i)} \left\vert \nabla \tilde{u}_n\right\vert^2 \rightarrow 0$ for  $i=1,2$ and  all $0<\alpha <1$, since there is no concentration to the boundary.

In particular, using again Proposition \ref{epsreg}, we have that for $i=1,2$,
\begin{equation} \label{estnei} \forall z\in \mathbb{D}_{\frac{\delta_n^i}{2}}(x_n^i)\cap \D_+ , \left\vert \tilde{u}_n(x_n^i) - \tilde{u}_n(z) \right\vert < \frac{\eta}{4} \hskip.1cm. \end{equation}

Now, we are going to extend $u_n$ to the whole plane. First we set for all $(x,y)\in \D_-$,

$$u_n(x,y)=-3u_n(x,-y)+4u_n\left(x,-\frac{y}{2}\right).$$
Then, the new $u_n \in W^{1,2}(\D,\R)$ and there exists $C>0$ independent of $n$, see section 5.4 of \cite{Evans} such that 
$$\Vert \nabla u_n\Vert_{L^2(\D)}\leq C \Vert \nabla u_n\Vert_{L^2(\D_+)} .$$
Then, we extend $\tilde{u}_n$ by an inversion on $\mathbb{R}^2$ so that $\tilde{u}_n(z) = \tilde{u}_n(\frac{1}{z})$ if $z\in \mathbb{R}^2\setminus\mathbb{D}$, which just doubles the energy of the map and keeps the energy small around the $x_n^i$ .\\

Then, for $i=1,2$, by the Courant-Lebesgue lemma, we take $ \frac{\delta_n^i}{4} < r_n^i < \frac{\delta_n^i}{2}$ such that
\begin{equation} \label{CLr1} \int_{\partial \mathbb{D}_{r_n^i}(x_n^i)} \left\vert \partial_\theta \tilde{u}_n\right\vert^2 d\theta \leq \frac{1}{\ln 2} \int_{\mathbb{D}_{\frac{\delta_n^i}{2}}(x_n^i)} \left\vert \nabla \tilde{u}_n\right\vert^2 
\end{equation}
\begin{equation} \label{CLr2}
\forall p,p' \in \partial \mathbb{D}_{r_n^i}(x_n^i),   \left\vert \tilde{u}_n(p) - \tilde{u}_n(p') \right\vert^2 \leq \frac{\pi}{\ln 2} \int_{\mathbb{D}_{\frac{\delta_n^i}{2}}(x_n^i)} \left\vert \nabla \tilde{u}_n\right\vert^2 
\end{equation}

Moreover, for $i=1,2$, again by the Courant-Lebesgue lemma, we take $ 2\delta_n^i < R_n^i < 4\delta_n^i$ such that
\begin{equation} \label{CLR1} \int_{\partial \mathbb{D}_{R_n^i}(x_n^i)} \left\vert \partial_\theta \tilde{u}_n\right\vert^2 d\theta \leq \frac{1}{\ln 2} \int_{\mathbb{D}_{4\delta_n^i}(x_n^i)} \left\vert \nabla \tilde{u}_n\right\vert^2 
\end{equation}
\begin{equation} \label{CLR2}
\forall q,q' \in \partial \mathbb{D}_{R_n^i}(x_n^i),   \left\vert \tilde{u}_n(q) - \tilde{u}_n(q') \right\vert^2 \leq \frac{\pi}{\ln 2} \int_{\mathbb{D}_{4\delta_n^i}(x_n^i)} \left\vert \nabla \tilde{u}_n\right\vert^2 
\end{equation}
Then, for $i=1,2$, we choose 
$$p_n^i \in \partial \mathbb{D}_{r_n^i}(x_n^i)\cap\D_+$$ 
and
$$q_n^i \in \partial \mathbb{D}_{R_n^i}(x_n^i) \cap A \hskip.1cm.$$
Those intersections are not empty since $r_n^i\leq \frac{\delta_i^n}{2}$ and $R_n^i \geq2\delta_n^i$.

For $i=1,2$, we let $v_n^i$ be the solution of 
\[
\left \{
\begin{array}{c  c c}
    \Delta v_n^i = 0 & \hbox{ in } \mathbb{D}_{R_n^i}(x_n^i)\setminus \mathbb{D}_{r_n^i}(x_n^i) \\
    v_n^i = \tilde{u}_n(p_n^i) & \hbox{ on } \partial \mathbb{D}_{r_n^i}(x_n^i) \\
    v_n^i = \tilde{u}_n(q_n^i) & \hbox{ on } \partial \mathbb{D}_{R_n^i}(x_n^i)
\end{array}
\right.
\]
and $w_n^i$ be the solution of
\[
\left \{
\begin{array}{c c c}
    \Delta w_n^i = 0 & \hbox{ in } \mathbb{D}_{R_n^i}(x_n^i)\setminus \mathbb{D}_{r_n^i}(x_n^i) \\
    w_n^i = \tilde{u}_n & \hbox{ on } \partial \mathbb{D}_{r_n^i}(x_n^i) \\
    w_n^i = \tilde{u}_n & \hbox{ on } \partial \mathbb{D}_{R_n^i}(x_n^i)
\end{array}
\right.
\]
Since $v_n^i - w_n^i$ is a harmonic map, and thanks to \eqref{CLr1}, \eqref{CLr2}, \eqref{CLR1}, \eqref{CLR2} and the fact the conformal class of the annuli is bounded since $4\leq R_n^i/r_n^i\leq 16 $, for $i=1,2$,
\begin{equation} \label{eqvw} \int_{\mathbb{D}_{R_n^i}(x_n^i)\setminus \mathbb{D}_{r_n^i}(x_n^i)} \left\vert\nabla v_n^i \right\vert^2 =  \int_{\mathbb{D}_{R_n^i}(x_n^i)\setminus \mathbb{D}_{r_n^i}(x_n^i)} \left\vert\nabla w_n^i \right\vert^2 + o(1) \end{equation}
But it is easy to notice that for $i=1,2$,
\begin{equation} \label{energyvni} \int_{\mathbb{D}_{R_n^i}(x_n^i)\setminus \mathbb{D}_{r_n^i}(x_n^i)} \left\vert\nabla v_n^i \right\vert^2 = \frac{2\pi}{\ln\left(\frac{R_n^i}{r_n^i}\right)} \left\vert \tilde{u}_n(p_n^i) -\tilde{u}_n(q_n^i) \right\vert^2 \hskip.1cm.\end{equation}
Then,
\begin{eqnarray*}
\eta & \leq & \left\vert \tilde{u}_n(x_n^1) -\tilde{u}_n(x_n^2) \right\vert \\
& \leq & \left( \sum_{i=1}^2  \left( \left\vert \tilde{u}_n(x_n^i) -\tilde{u}_n(p_n^i) \right\vert + \left\vert \tilde{u}_n(p_n^i) -\tilde{u}_n(q_n^i) \right\vert \right) \right)  + \left\vert \tilde{u}_n(q_n^1) -\tilde{u}_n(q_n^2) \right\vert \\
& \leq & \frac{\eta}{2} + \left(\frac{\ln\left(\frac{R_n^i}{r_n^i}\right)}{2\pi}\right)^{\frac{1}{2}} \left(  \sum_{i=1}^2 \left(\int_{\mathbb{D}_{R_n^i}(x_n^i)\setminus \mathbb{D}_{r_n^i}(x_n^i)} \left\vert\nabla v_n^i \right\vert^2 \right)^{\frac{1}{2}}  \right) + o(1) \\
& \leq & \frac{\eta}{2} + \left(\frac{\ln\left(\frac{R_n^i}{r_n^i} \right)}{2\pi}\right)^{\frac{1}{2}} \left(  \sum_{i=1}^2 \left(\int_{\mathbb{D}_{R_n^i}(x_n^i)\setminus \mathbb{D}_{r_n^i}(x_n^i)} \left\vert\nabla w_n^i \right\vert^2 \right)^{\frac{1}{2}} + o(1) \right) +o(1) \\
& \leq & \frac{\eta}{2} + 2\left(\frac{\ln\left(\frac{R_n^i}{r_n^i} \right)}{2\pi}\right)^{\frac{1}{2}} \left(\int_{\mathbb{D}_{4\delta_n^i}(x_n^1) \cup \mathbb{D}_{4\delta_n^i}(x_n^2)   } \left\vert\nabla \tilde{u}_n \right\vert^2 \right)^{\frac{1}{2}} + o(1) \\
& \leq &  \frac{\eta}{2} + o(1)
\end{eqnarray*}
where we used 
\begin{itemize}
\item in the third inequality that \eqref{estnei} and the definition of $p_n^i$, \eqref{energyvni} and that $\left\vert q_n^1 - q_n^2 \right\vert \leq R_n^1 + R_n^2 + \frac{1}{n} \to 0$ as $n\to + \infty$ with $q_n^i \in A$ so that $\left\vert \tilde{u}_n(q_n^1) -\tilde{u}_n(q_n^2) \right\vert = \left\vert u_n(q_n^1) - u_n(q_n^2) \right\vert = o(1)$,
\item in the fourth inequality \eqref{eqvw} 
\item in the fifth inequality the definition of $w_n^i$ 
\item in the last inequality \textit{STEP 1} and $4\leq R_n^i/r_n^i\leq 16 $.
\end{itemize}
This leads to a contradiction and achieves the proof of STEP3.

\vspace{4mm}

Gathering \textit{STEP 2} and \textit{STEP 3}, the proof of Proposition \ref{continuity} is complete.

\hfill $\diamondsuit$
\subsection{The replacement procedure}
Let $\mathcal{E}$ the set of finite families of disjoint element of $\mathfrak{B}$ and $\mathfrak{HB}$ where 
$$\mathfrak{B}=\{ \overline{B(a,r)}\hbox{ s.t. } B(a,r)\subset \mathbb{D}\}$$
and 
$$\mathfrak{HB}=\{\overline{B(a,r)\cap \mathbb{D}} \hbox{ s.t. } \partial\mathbb{D} \hbox{ and } \partial B(a,r) \hbox{ intersect orthogonally}\}.$$

For any element $B$ of $ \mathfrak{HB}$ we denote $A_B=\partial B\cap \mathbb{D}$.

For $\alpha\in (0,1]$, and $\mathcal{B}\in \mathcal{E}$ we denote by $\alpha\mathcal{B}$ the collection of disks and concentric to those of $\mathcal{B}$ with radius dilated by $\alpha$ and half disks which meet the boundary orthogonally, with center collinear to those of $\mathcal{B}$ and with radius dilated by $\alpha$..

For $u:\mathbb{D}\rightarrow N$ with $u(\partial{\mathbb{D}}) \subset M$ and $\mathcal{B} \in \mathcal{E}$ such that the energy of $u$ on $\bigcup_{B\in \mathcal{B}} B$ is less than the  $\eps_0$ of theorem \ref{econvexb}, then we denote by $H(u,\mathcal{B}):\mathbb{D}\rightarrow N$ the map that coincides 
\begin{itemize}
\item on $\mathbb{D}\setminus \bigcup_{B\in \mathcal{B}} B$ with $u$
\item on $\bigcup_{B\in \mathcal{B}} B$ with the (unique) energy minimizing map from $\bigcup_{B\in \mathcal{B}} B$ to $N$ that agrees with $u$ on $\bigcup_{B\in \mathcal{B}} \partial B \setminus \partial{\mathbb{D}}$ with the constraint that it lies in $M$ on $\bigcup_{B\in \mathcal{B}} \partial B \cap \partial{\mathbb{D}}$.
\end{itemize}

We will also denote by induction 
$$H(u,\mathcal{B}_1,\mathcal{B}_2,\cdots,\mathcal{B}_{\tau}) = H(H(u,\mathcal{B}_1,\mathcal{B}_2,\cdots,\mathcal{B}_{\tau-1}),\mathcal{B}_{\tau}) \hskip.1cm.$$

We will need the following proposition :

\begin{prop} \label{exchangeharm}
There is a constant $\kappa >0$ depending on $M$ and $N$ such that if $u\in \mathcal{C}^0(\overline{\mathbb{D},N}) \cap W^{1,2}(\mathbb{D},N)$ with $u(\partial{\mathbb{D}}) \subset M$ and if $\mathcal{B}_1$ and $\mathcal{B}_2$ lie in $\mathcal{E}$ so that the energy of $u$ on $\bigcup_{B\in \mathcal{B}} B$ is less than $\frac{\eps_0}{3}$, then
\begin{equation} \label{exchange1} E(u)-E(H(u, \frac{1}{2}\mathcal{B}_2)) \leq E(H(u,\mathcal{B}_1))-E(H(u, \mathcal{B}_1,\mathcal{B}_2)) +\frac{(E(u)-E(H(u, \mathcal{B}_1)))^{\frac{1}{2}}}{\kappa} \end{equation}
and
\begin{equation} \label{exchange2} E(H(u, \mathcal{B}_1))-E(H(u, \mathcal{B}_1, \frac{1}{2}\mathcal{B}_2)) \leq E(u)-E(H(u, \mathcal{B}_2)) +\frac{(E(u)-E(H(u, \mathcal{B}_1)))^{\frac{1}{2}}}{\kappa} \hskip.1cm. \end{equation}
\end{prop}

In order to prove Proposition \ref{exchangeharm}, we use the following lemmas. The first one was proved by Colding and Minicozzi, see lemma 3.11 \cite{CM2008}, and the second one is an adaptation in the case of half-disks. Again, a careful attention is given to the setting of competitors which are admissible functions (in particular, sending $I$ into $M$).

\begin{lem} \label{extensionestimates} There is $\eta >0 $ and a large constant $K$ depending on $N$ such that for any $R>0$ and for any $f,g \in W^{1,2}(\partial \mathbb{D}_R, N)$, if $f$ and $g$ agree at one point, $f\neq g$ and\footnote{$\nabla_{\theta}f$ means $\frac{df}{rd\theta}$.}
$$ R \int_{\partial \mathbb{D}_R} \left\vert \nabla_{\theta}(f - g) \right\vert^2 \,d\sigma \leq \eta^2,$$
then we can find $\rho \in \left(0,\frac{R}{2}\right]$ and a map $w \in \mathcal{C}^0\cap W^{1,2}(\mathbb{D}_R \setminus \mathbb{D}_{R-\rho}, N)$ with $w_{\vert \partial \mathbb{D}_{R-\rho}} = f\left(\frac{R}{R-\rho} .\right)$ and $w_{\vert \partial \mathbb{D}_{R}} = g$ which satisfies the estimates
\begin{equation*}\int_{\mathbb{D}_R\setminus \mathbb{D}_{R-\rho}} \left\vert\nabla w \right\vert^2 \leq K \left( R \int_{\partial \mathbb{D}_R} \left\vert \nabla_\theta (f - g) \right\vert^2 \,d\sigma\right)^{\frac{1}{2}} \left( R \int_{\partial \mathbb{D}_R} \left( \left\vert \nabla_\theta f \right\vert^2 +  \left\vert \nabla_\theta g \right\vert^2\,d\sigma \right) \right)^{\frac{1}{2}}.
\end{equation*}
\end{lem}

\begin{lem} \label{bextensionestimates} There is $\eta >0 $ and a large constant $K$ depending on $N$ and $M$ such that for any $R>0$ and for any $f,g \in W^{1,2}(A_R, N)$ with $f(\partial A_R) \subset M$ and $g(\partial A_R) \subset M$, where $A_R=\partial B(0,R)\cap\{y\geq 0 \}$, if $f$ and $g$ agree at one point, $f\neq g$ and  
$$ R \int_{A_R} \left\vert \nabla_{\theta}(f - g) \right\vert^2 \,d\sigma \leq \eta^2 ,$$
where $\nabla_{\theta}f$ means $\frac{df}{rd\theta}$,then we can find $\rho \in (0,\frac{R}{2}]$ and a map $w \in \mathcal{C}^0\cap W^{1,2}(\mathbb{D}_R^+ \setminus \mathbb{D}_{R-\rho}^+, N)$ with $w(I_R\setminus I_{R-\rho}) \subset M$, $w_{\vert A_{R-\rho}} = f\left(\frac{R}{R-\rho} .\right)$ and $w_{\vert A_{R}} = g$ which satisfies the estimates
\begin{equation*}
\int_{\mathbb{D}_R^+\setminus \mathbb{D}_{R-\rho}^+} \left\vert\nabla w \right\vert^2 \leq K \left( R \int_{A_R} \left\vert \nabla_\theta (f - g) \right\vert^2 \,d\sigma \right)^{\frac{1}{2}} \left( R \int_{A_R} \left( \left\vert \nabla_\theta f \right\vert^2 +  \left\vert \nabla_\theta g \right\vert^2 \,d\sigma\right) \right)^{\frac{1}{2}}.
\end{equation*}
\end{lem}

{\it Proof of lemma \ref{bextensionestimates} :}\\
Since the statement is scale invariant, we can assume that $R=1$. We set
$$ \rho^2 = \frac{\int_{A_R} \left\vert \nabla_\theta(f-g) \right\vert^2}{8 \int_{A_R} \left( \left\vert \nabla_\theta f \right\vert^2 +  \left\vert \nabla_\theta g \right\vert^2 \right)} \hskip.1cm.$$
Since $f\neq g$ we have $0<\rho \leq \frac{1}{2}$. 

We define
\begin{equation} \label{defwhat} \hat{w}(r,\theta) = f(\theta) + \frac{r+\rho-1}{\rho}\left( g(\theta) - f(\theta) \right) \end{equation}
for $1-\rho \leq r \leq 1$. 

We have that
$$ dist(\hat{w}(r,0),M) \leq \left\vert \hat{w}(r,0) - f(0) \right\vert = \frac{r+\rho-1}{\rho} \left\vert g(0) - f(0) \right\vert \leq \left\vert g(0) - f(0) \right\vert \leq \sqrt{\pi}\eta$$
and similarly
$$ dist(\hat{w}(r,\pi),M) \leq \sqrt{\pi}\eta.$$

We choose $\eta \leq \frac{\delta}{\sqrt{\pi}}$ so that we can set :
\begin{equation} \label{defwtilde} \tilde{w}(r,\theta) = \hat{w}(r,\theta) + \frac{\pi-\theta}{\pi}\left( \pi_M(\hat{w}(r,0)) - \hat{w}(r,0)  \right) + \frac{\theta}{\pi}\left( \pi_M(\hat{w}(r,\pi)) - \hat{w}(r,\pi)  \right) \end{equation}
so that 
$$ \left\vert  \nabla\tilde{w} \right\vert^2 \leq 3 \left( \left\vert  \nabla\hat{w} \right\vert^2 + \left\vert  \nabla\left\{ \frac{\pi-\theta}{\pi}(\pi_M - id)( \hat{w}(r,0) )\right\}\right\vert^2  + \left\vert   \nabla \left\{ \frac{\theta}{\pi}(\pi_M - id)(\hat{w}(r,\pi))  \right\} \right\vert^2 \right) $$
We have that
\begin{eqnarray*}
\left\vert  \nabla\hat{w} \right\vert^2(r,\theta) & = & \left\vert \frac{g(\theta)-f(\theta)}{\rho} \right\vert^2 + \frac{1}{r^2} \left\vert \frac{1-r}{\rho}\nabla_\theta f(\theta) + \frac{r+\rho-1}{\rho}\nabla_\theta g(\theta) \right\vert^2 \\
& \leq & \frac{\pi}{\rho^2}\int_A \left\vert \nabla_\theta (f-g) \right\vert^2 + \frac{2}{r^2}\left( \left\vert \nabla_\theta f \right\vert^2 + \left\vert \nabla_\theta g \right\vert^2 \right) \hskip.1cm.
\end{eqnarray*} 
We also have that
\begin{eqnarray*}
\left\vert  \nabla\left\{ \frac{\pi-\theta}{\pi}(\pi_M - id)( \hat{w}(r,0) )\right\}\right\vert^2 & = & \left\vert  \partial_r\left\{ \frac{\pi-\theta}{\pi}(\pi_M - id)( \hat{w}(r,0) )\right\}\right\vert^2 \\
& & + \frac{1}{r^2 \pi^2} \left\vert (\pi_M - id)( \hat{w}(r,0) )\right\vert^2
\end{eqnarray*}
with
\begin{eqnarray*}
\left\vert (\pi_M - id)( \hat{w}(r,0) )\right\vert^2 & = & \left\vert \int_{1-\rho}^r  \partial_s \left\{ (\pi_M - id)( \hat{w}(s,0) )\right\}ds \right\vert^2 \\
& \leq & \left\| d\pi_M - id \right\|_{L^{\infty}(M_{\hat{\delta}})}^2 \left( \int_{1-\rho}^r  \left\vert\partial_s  \hat{w}(s,0) \right\vert ds \right)^2  \\
& \leq & 9 \left\vert f(0)-g(0) \right\vert^2 \\
& \leq & 9\pi \int_A \left\vert \nabla_\theta (f- g) \right\vert^2 \\
& \leq & 18\pi \int_A \left(\left\vert \nabla_\theta f \right\vert^2 + \left\vert \nabla_\theta g \right\vert^2\right)
\end{eqnarray*}
and 
\begin{eqnarray*}
\left\vert  \partial_r\left\{ \frac{\pi-\theta}{\pi}(\pi_M - id)( \hat{w}(r,0) )\right\}\right\vert^2 &\leq & \left\| d\pi_M - id \right\|_{L^{\infty}(M_{\hat{\delta}})}^2 \left\vert  \partial_r \hat{w}(r,0) \right\vert^2 \\
& \leq & \frac{9\pi}{\rho^2} \int_A \left\vert \nabla_\theta (f-g) \right\vert^2
\end{eqnarray*}
so that
$$
\left\vert  \nabla\left\{ \frac{\pi-\theta}{\pi}(\pi_M - id)( \hat{w}(r,0) )\right\}\right\vert^2  \leq  C_0 \left( \frac{1}{\rho^2} \int_A \left\vert \nabla_\theta (f-g) \right\vert^2 + \frac{1}{r^2} \int_A \left(\left\vert \nabla_\theta f \right\vert^2 + \left\vert \nabla_\theta g \right\vert^2\right) \right)
$$
for some constant $C_0$. The same computation gives
$$
\left\vert  \nabla\left\{ \frac{\theta}{\pi}(\pi_M - id)( \hat{w}(r,\pi) )\right\}\right\vert^2  \leq  C_0 \left( \frac{1}{\rho^2} \int_A \left\vert \nabla_\theta (f- g) \right\vert^2 + \frac{1}{r^2} \int_A \left(\left\vert \nabla_\theta f \right\vert^2 + \left\vert \nabla_\theta g \right\vert^2\right) \right)
$$
so that gathering all the previous inequalities, we get a constant $C_1>0$ such that
$$ \int_0^{\pi} \left\vert  \nabla\tilde{w} \right\vert^2(r,\theta)d\theta \leq C_1  \left( \frac{1}{\rho^2} \int_A \left\vert \nabla_\theta (f- g) \right\vert^2 + \frac{1}{r^2} \int_A \left(\left\vert \nabla_\theta f \right\vert^2 + \left\vert \nabla_\theta g \right\vert^2\right) \right) $$
for any $1-\rho\leq r \leq 1$. Now,
\begin{eqnarray*}
\int_{\mathbb{D}_+ \setminus \mathbb{D}_{1-\rho}} \left\vert  \nabla\tilde{w} \right\vert^2 & \leq & C_1 \left(\int_{1-\rho}^1 \left( \frac{r}{\rho^2} \int_A \left\vert \nabla_\theta (f- g) \right\vert^2 + \frac{1}{r} \int_A \left(\left\vert \nabla_\theta f \right\vert^2 + \left\vert \nabla_\theta g \right\vert^2\right) \right)\right) \\
& \leq & C_1 \left( \frac{1}{\rho} \int_A \left\vert \nabla_\theta (f- g) \right\vert^2 + \rho \int_A \left(\left\vert \nabla_\theta f \right\vert^2 + \left\vert \nabla_\theta g \right\vert^2\right) \right) \\
& \leq & C_1  \left(  \int_{A} \left\vert \nabla_\theta(f-g) \right\vert^2 \right)^{\frac{1}{2}} \left(  \int_{A} \left( \left\vert \nabla_\theta f \right\vert^2 +  \left\vert \nabla_\theta g \right\vert^2 \right) \right)^{\frac{1}{2}} \hskip.1cm.
\end{eqnarray*}

Thanks to the definition of $\tilde{w}$ (see \eqref{defwtilde} and \eqref{defwhat}) and the assumptions $f(0),g(0),f(\pi),g(\pi) \in M$ we have that $\tilde{w}(I\setminus I_{1-\rho}) \subset M$, $\tilde{w}_{\vert A_{1-\rho}} = f((1-\rho)^{-1} .)$ and $\tilde{w}_{\vert A} = g$. Moreover, letting $1-\rho \leq r \leq 1$ and $0 \leq\theta \leq \pi$, we have that
\begin{eqnarray*}
dist(\tilde{w}(r,\theta),N) & \leq & dist(\hat{w}(r,\theta),N) + dist(\hat{w}(r,0),M) + dist(\hat{w}(r,\pi),M) \\  
&\leq & \left\vert \hat{w}(r,\theta) - f(\theta) \right\vert + 2 \sqrt{\pi} \eta \\
& \leq & \frac{\rho-1+r}{\rho}\left\vert g(\theta) - f(\theta) \right\vert \\
& \leq & 3\sqrt{\pi} \eta
\end{eqnarray*}
So that letting $\eta \leq \frac{\delta}{3\sqrt{\pi}} $ we can set $ w = \pi_N (\tilde{w}) $. The projection on $N$ gives that $w(I\setminus I_{1-\rho}) \subset M$, $w_{\vert A_{1-\rho}} = f((1-\rho) .)$ and $w_{\vert A} = g$ and we obtain the inequality
\begin{eqnarray*}
\int_{\mathbb{D}_+ \setminus \mathbb{D}_{1-\rho}} \left\vert  \nabla w \right\vert^2 & \leq & 4 \int_{\mathbb{D}_+ \setminus \mathbb{D}_{1-\rho}} \left\vert  \nabla\tilde{w} \right\vert^2 \\
& \leq & 4 C_1  \left(  \int_{A} \left\vert \nabla_\theta(f-g) \right\vert^2 \right)^{\frac{1}{2}} \left(  \int_{A} \left( \left\vert \nabla_\theta f \right\vert^2 +  \left\vert \nabla_\theta g \right\vert^2 \right) \right)^{\frac{1}{2}}
\end{eqnarray*}  
which ends the proof of proposition \ref{bextensionestimates}.

\hfill $\diamondsuit$

\medskip

The following proof is analogous to the one given by Colding Minicozzi \cite{CM2008} (lemma 3.8) even if we slightly change the presentation and give a proof which works in the free-boundary case.

\medskip

{\it Proof of proposition \ref{exchangeharm} :}\\

Notice that the energy of $u$ on $\mathcal{B}_1 \cup \mathcal{B}_2$ is less than $\frac{2\eps_0}{3}$ so that the energy of $H(u,\mathcal{B}_1)$ on $\mathcal{B}_2$ is less than $\frac{2\eps_0}{3}$. For all the proof, we set
$$ \mathcal{B}_2^- = \{ B \in \mathcal{B}_2 ; \forall \tilde{B} \in \mathcal{B}_1 , \frac{1}{2}B \nsubseteq \tilde{B} \} $$
$$ \mathcal{B}_2^+ = \mathcal{B}_2 \setminus \mathcal{B}_2^- $$
We set $u_1 = H(u,\mathcal{B}_1)$ and we only prove \eqref{exchange1} since we can follow the same proof for \eqref{exchange2} switching the role of $u_1$ and $u$.

We can assume that $ 9 \int_{\mathbb{D}} \left\vert \nabla u_1 - \nabla u \right\vert^2 \leq \eta^2$ where $\eta$ is given by lemmas \ref{extensionestimates} and \ref{bextensionestimates}. If not, we use the convexity inequalities in Theorems \ref{econvexb} and \ref{econvex} to get
$$ \int_{\mathbb{D}} \left\vert \nabla u \right\vert^2 - \int_{\mathbb{D}} \left\vert \nabla u_1  \right\vert^2 \geq \frac{\eta^2}{18} = \kappa^2 \eps_0^2 $$
setting $\kappa = \frac{\eta}{18^{\frac{1}{2}}\eps_0}$ so that \eqref{exchange1} is true.

\vspace{4mm}
{\it STEP 1 :} Let $B \in \mathcal{B}_2^-$. We prove that there is a universal constant $C>0$ such that
\begin{eqnarray*}
\int_B \left\vert \nabla u_1 \right\vert^2 - \int_B \left\vert \nabla H(u_1,B) \right\vert^2 & \geq & \int_{\frac{1}{2}B} \left\vert \nabla u \right\vert^2 - \int_{\frac{1}{2}B} \left\vert \nabla H(u,\frac{1}{2}B) \right\vert^2 \\
& & - C \left( \int_B \left\vert\nabla u\right\vert^2+\left\vert\nabla u_1\right\vert^2\right)^{\frac{1}{2}} \left( \int_B \left\vert \nabla(u-u_1)\right\vert^2 \right)^{\frac{1}{2}}.
\end{eqnarray*}
 
We denote $R$ the radius of $B$. Then there is $r \in\left[\frac{3R}{4},R\right]$  such that
$$ \int_{\partial B_r\cap \D} \left\vert \nabla u_1-\nabla u \right\vert^2 \leq \frac{9}{R} \int_{\frac{3R}{4}}^{R} \int_{\partial B_s\cap \D} \left\vert \nabla u_1-\nabla u \right\vert^2 ds \leq  \frac{9}{r} \int_{B_R} \left\vert \nabla u_1-\nabla u \right\vert^2 $$
 and
 $$ \int_{\partial B_r\cap \D} \left\vert \nabla u \right\vert^2 + \left\vert \nabla u_1 \right\vert^2 \leq \frac{9}{R} \int_{\frac{3R}{4}}^{R} \int_{\partial B_s\cap \D} \left\vert  \nabla u \right\vert^2 + \left\vert \nabla u_1 \right\vert^2 ds \leq  \frac{9}{r} \int_{B_R} \left\vert \nabla u \right\vert^2 + \left\vert \nabla u_1 \right\vert^2 $$
Since $B\in \mathcal{B}_2^-$ and $r>\frac{R}{2}$, then $\partial B_r$ contains a point outside every ball of $\mathcal{B}_1$ : $u$ and $u_1$ coincide at such a point. By lemmas \ref{bextensionestimates} and \ref{extensionestimates}, if $u_1 \neq u$ on $\partial B_r$, there is $\rho \in (0,\frac{R}{2}]$ and $w : B_r\setminus B_{r-\rho} \rightarrow N$ such that $w(r,\theta) = u_1(r,\theta)$, $w(r-\rho,\theta)=u(r,\theta)$ and
\begin{equation} \label{ineqw}\int_{B_r\setminus B_{R-\rho}} \left\vert\nabla w\right\vert^2 \leq K \left(\int_{B_R} \left\vert \nabla u_1-\nabla u \right\vert^2\right)^{\frac{1}{2}}\left(\int_{B_R}\left\vert \nabla u \right\vert^2 + \left\vert \nabla u_1 \right\vert^2\right)^{\frac{1}{2}} \end{equation}
and if $u=u_1$ on $\partial B_r$, we let $\rho = 0$. In fact the $w$ furnished by lemma \ref{extensionestimates} is defined on $\D_r⁺\setminus \D_{r-\rho}⁺$ but can easily push it to $B_r\setminus B_{r-\rho}$ through the conformal diffeomorphism from the half plane to the disk, which leave the considered norm unchanged. Then we set
$$ v(x) = \begin{cases}
u_1(x) & \text{ if }x\in B_R\setminus B_r\\
w(x) & \text{ if }x\in B_r\setminus B_{r-\rho} \\
H(u,B_r)\left(\frac{rx}{r-\rho}\right) & \text{ if }x\in  B_{r-\rho}
\end{cases} $$
Then $v \in W^{1,2}(B)$ and
\begin{eqnarray*}
\int_{B_R}\left\vert \nabla H(u_1,B_R)\right\vert^2 & \leq & \int_{B_R}\left\vert \nabla v\right\vert^2 \\
&=& \int_{B_R\setminus B_r}\left\vert \nabla u_1 \right\vert^2 + \int_{B_r\setminus B_{r-\rho} }\left\vert \nabla w \right\vert^2 + \int_{B_{r-\rho} }\left\vert \nabla H(u,B_r)\left(\frac{rx}{r-\rho}\right) \right\vert^2 \\
&=& \int_{B_R\setminus B_r}\left\vert \nabla u_1 \right\vert^2 + \int_{B_r\setminus B_{r-\rho} }\left\vert \nabla w \right\vert^2 + \int_{B_{r} }\left\vert \nabla H(u,B_r) \right\vert^2 \\
\end{eqnarray*}
Then,
\begin{eqnarray*}
\int_{B_R}\left\vert \nabla u_1\right\vert^2 - \int_{B_R}\left\vert \nabla H(u_1,B_R)\right\vert^2  & \geq & \int_{B_r}\left\vert \nabla u_1\right\vert^2 - \int_{B_{r} }\left\vert \nabla H(u,B_r) \right\vert^2 - \int_{B_r\setminus B_{r-\rho} }\left\vert \nabla w \right\vert^2  \\
 & \geq & \int_{B_r}\left\vert \nabla u\right\vert^2 - \int_{B_{r} }\left\vert \nabla H(u,B_r) \right\vert^2 \\
 & & - (K+\sqrt{2})\left(\int_{B_R} \left\vert \nabla u_1-\nabla u \right\vert^2\right)^{\frac{1}{2}}\left(\int_{B_R}\left\vert \nabla u \right\vert^2 + \left\vert \nabla u_1 \right\vert^2\right)^{\frac{1}{2}} \hskip.1cm.
\end{eqnarray*}
Here we used \eqref{ineqw} and that
\begin{eqnarray*}
\int_{B_r}\left\vert \nabla u_1\right\vert^2 - \int_{B_r}\left\vert \nabla u\right\vert^2 &\geq & - \int_{B_r} \left( \left\vert\nabla u\right\vert+ \left\vert\nabla u_1\right\vert\right)\left\vert\nabla u-\nabla u_1\right\vert \\
& \geq & -\sqrt{2}\left(\int_{B_r} \left\vert \nabla u_1-\nabla u \right\vert^2\right)^{\frac{1}{2}}\left( \int_{B_r}\left\vert \nabla u \right\vert^2 + \left\vert \nabla u_1 \right\vert^2\right)^{\frac{1}{2}} \hskip.1cm.
\end{eqnarray*}
We get STEP 1 with $C = K+\sqrt{2}$ noticing that
$$ \int_{B_{r} }\left\vert \nabla H(u,B_r) \right\vert^2 \leq \int_{B_{r} \setminus  B_{\frac{R}{2} }}\left\vert \nabla u\right\vert^2 + \int_{B_{\frac{R}{2} }}\left\vert \nabla H(u,B_{\frac{R}{2}}) \right\vert^2 \hskip.1cm.$$
\vspace{4mm}

{\it STEP 2 :} We have
\begin{equation} \label{inequalitystep2}  E(u_1)- E(H(u_1,\mathcal{B}_2^-)) \geq E(u)-E(H(u, \frac{1}{2}\mathcal{B}_2^-)) - \frac{(E(u)-E(u_1))^{\frac{1}{2}}}{\kappa} \end{equation}

Indeed, we sum the inequality of STEP 1 for all $B \in \mathcal{B}_2^-$, and we set $\kappa = \frac{\eps_0^{\frac{1}{2}}}{C}$. We get this inequality using that
$$ \sum a_j b_j \leq \left(\sum a_j^2\right)^{\frac{1}{2}}\left(\sum b_j^2\right)^{\frac{1}{2}} $$
and that by the convexity inequality (Theorem \ref{econvexb}and \ref{econvex}),
$$ \int_{\mathbb{D}} \left\vert \nabla (u-u_1)\right\vert^2 \leq \frac{1}{2} \left( E(u) - E(u_1) \right) $$
\vspace{4mm}

{\it STEP 3 :} We conclude by the proof of \eqref{exchange1}. We apply 
$$ E(H(u,\mathcal{B}_1,\mathcal{B}_2^+)) \leq  E(H(u,\mathcal{B}_1,\frac{1}{2}\mathcal{B}_2^+)) = E(H(u,\mathcal{B}_1)) \leq E(H(u, \frac{1}{2}\mathcal{B}_2^+)) \hskip.1cm, $$
to obtain
\begin{equation} \label{ineqforstep3} E(u) - E\left(H\left(u,\frac{1}{2}\mathcal{B}_2^+\right)\right) \leq E(u) - E(u_1) + E(u_1) - E\left(H\left(u_1,\mathcal{B}_2^+\right)\right) \hskip.1cm.\end{equation}
We then apply \eqref{inequalitystep2} and \eqref{ineqforstep3} on the right side of the equality
$$ E(u) - E\left(H\left(u,\frac{1}{2}\mathcal{B}_2\right)\right) =  E(u) - E\left(H\left(u,\frac{1}{2}\mathcal{B}_2^-\right)\right) + E(u) - E\left(H\left(u,\frac{1}{2}\mathcal{B}_2^+\right)\right) $$
and we get \eqref{exchange1} with
$$ E(u_1) - E\left(H\left(u_1,\mathcal{B}_2\right)\right) =E(u_1) - E\left(H\left(u_1,\mathcal{B}_2^+\right)\right) + E(u_1) - E\left(H\left(u_1,\mathcal{B}_2^-\right)\right) \hskip.1cm.$$

\hfill $\diamondsuit$

Then, by Besicovitch covering theorem, we know that for $k\in\mathbb{N}$ there exists $\tau$, depending only on $k$, such that for any covering of a compact subset $K$ of $\mathbb{R}^{k-2}$ by a family balls $\{B(x,r_x)\}_{x\in K}$ with $r_x>0$, then there is a finite sub-cover $\{B_i\}_{i=1\cdots m}$ such that
$$1 \leq \sum_{i=1}^m \mathbf{1}_{B_i} \leq \tau \hskip.1cm.$$

Let $\sigma : \mathbb{B}^{k-2} \rightarrow \mathcal{A}$ be a sweepout where we denote
$$ \mathcal{A} = \{u\in W^{1,2}\cap\mathcal{C}^0(\overline{\mathbb{D}},N); u(\partial\mathbb{D}) \subset M\} $$
the set of admissible maps. We set for some $t\in \mathbb{B}^k$

$$ e_{\sigma,\eps}(t) = \sup \left\{  E(\sigma_t) - E(H(\sigma_t, \frac{1}{2^{\tau-1}}\mathcal{B}))  ; \mathcal{B}\in \mathcal{E} , \int_{\bigcup_{B\in\mathcal{B}} B} \left\vert \nabla \sigma_t \right\vert^2 \leq \eps \right\}.$$

\begin{prop} \label{ballinBn} Let $t\in \mathbb{B}^{k-2}$. If $\sigma_t$ is not a harmonic map with free boundary and $0<\eps\leq\eps_0$, then there is a ball $C_t$ in $\mathbb{B}^{k-2}$ centered in $t$ such that
$$ \forall s \in 2C_t, e_{\sigma,\frac{\eps}{2}}(s) \leq 2 e_{\sigma,\eps}(t)  \hskip.1cm. $$
\end{prop}

{\it Proof of proposition \ref{ballinBn} :}\\

Notice that since $\sigma_t$ is not harmonic, $e_{\sigma,\eps}(t)>0$. By proposition \ref{continuity}, and since $t \in \mathbb{B}^k \mapsto \sigma_t \in \mathcal{A}$ is continuous, then $t\mapsto E(H(\sigma_t,\frac{1}{2^{\tau-1}}\mathcal{B}))$ is also continuous. Let a ball $C_t$ centered in $t$ such that for any $s \in 2C_t$
\begin{equation} \label{continuousHenergy}
\left\vert E(H(\sigma_t,\frac{1}{2^{\tau-1}}\mathcal{B})) - E(H(\sigma_s,\frac{1}{2^{\tau-1}}\mathcal{B})) \right\vert \leq \frac{e_{\sigma,\eps}(t)}{2}
\end{equation}
and
\begin{equation} \label{continuousenergy}
\int_{\mathbb{D}} \left\vert \left\vert \nabla\sigma_t\right\vert^2 -  \left\vert\nabla\sigma_s\right\vert^2\right\vert \leq \min\left\{\frac{\eps}{2},\frac{e_{\sigma,\eps}(t)}{2}\right\}\hskip.1cm. 
\end{equation}
Now, let $s\in 2 C_t$ and $\mathcal{B}\in\mathcal{E}$ such that $\int_{\bigcup_{B\in\mathcal{B}}B}\left\vert\nabla\sigma_s\right\vert^2 \leq \frac{\eps}{2}$. Then, by  \eqref{continuousenergy},
\begin{equation} \label{condeps}\int_{\bigcup_{B\in\mathcal{B}}B}\left\vert\nabla\sigma_t\right\vert^2 \leq \eps \end{equation}
and by \eqref{continuousHenergy} and \eqref{continuousenergy},
\begin{equation} \label{diffenergyst} \left\vert E(\sigma_s) - E(H(\sigma_s,\frac{1}{2^{\tau-1}}\mathcal{B})) - ( E(\sigma_t) - E(H(\sigma_t,\frac{1}{2^{\tau-1}}\mathcal{B})) ) \right\vert \leq e_{\sigma,\eps}(t) \hskip.1cm. \end{equation}
With \eqref{condeps} and \eqref{diffenergyst} and the definition of $e_{\sigma,\eps}(t)$, we get that 
$$ \left\vert E(\sigma_s) - E(H(\sigma_s,\frac{1}{2^{\tau-1}}\mathcal{B})) \right\vert \leq 2 e_{\sigma,\eps}(t)$$
and we take the supremum on $\mathcal{B}$.
\hfill $\diamondsuit$

\begin{prop} \label{deffamilies} We assume that $W>0$. Let $\tilde{\sigma}$ be a sweepout such that for any $t \in \mathbb{B}^{k-2}$, if $\tilde{\sigma}_t$ is a harmonic map with free boundary, then $\tilde{\sigma}_t$ is constant. Then, there are families $\mathcal{B}_1,\cdots,\mathcal{B}_m \in \mathcal{E}$ and continuous functions $r_1,\cdots,r_m : \mathbb{B}^{k-2} \rightarrow [0,1]$ such that for any $t\in \mathbb{B}^{k-2}$,
\begin{itemize}
\item $\sharp \{j ; r_j(t) > 0 \} \leq \tau $ and for any such $j$, $\int_{\bigcup_{B\in\mathcal{B}_j} B} \left\vert \nabla \tilde{\sigma}_t \right\vert^2 \leq \frac{\eps_0}{3^{\tau-1}}$
\item If $E(\tilde{\sigma}_t) \geq \frac{W}{2}$, then there is $j(t)$ such that
$$ E(\tilde{\sigma}_t) - E\left(H\left(\tilde{\sigma}_t,\frac{r_{j(t)}(t)}{2^{\tau-1}}B_{j(t)}\right)\right) \geq \frac{e_{\tilde{\sigma}, \frac{\eps_0}{3^{\tau+1}}}(t)}{8}$$
\end{itemize}

\end{prop}

{\it Proof of proposition \ref{deffamilies} :}\\

We set $A = \{t\in \mathbb{B}^{k-2} ; E(\tilde{\sigma}_t)\geq \frac{W}{2}\}$. Let $t\in A$. Then, there is $\mathcal{B}_t \in \mathcal{E}$ such that
\begin{itemize}
\item $\int_{\bigcup_{B\in\mathcal{B}_t} B} \left\vert \nabla \tilde{\sigma}_t \right\vert^2 \leq \frac{\eps_0}{3^{\tau}}$,
\item $E(\tilde{\sigma}_t) - E\left(H\left(\tilde{\sigma}_t,\frac{1}{2^{\tau-1}}B_{t}\right)\right) \geq \frac{e_{\tilde{\sigma}, \frac{\eps_0}{3^{\tau}}}(t)}{2}$
\end{itemize}
By Proposition \ref{ballinBn}, there is a closed ball $C_t$ centered at $t$ such that for any $s\in 2C_t$,
$$ e_{\sigma,\frac{\eps_0}{3^{\tau+1}}}(s) \leq 2 e_{\sigma,\frac{\eps_0}{3^{\tau}}}(t)  \hskip.1cm. $$
By continuity of $s\mapsto \tilde{\sigma}_s$ in $\mathcal{C}^0\cap W^{1,2}$ we reduce $C_t$ so that
\begin{itemize}
\item $\tilde{\sigma}_s$ has energy less than $\frac{\eps_0}{3^{\tau-1}}$ on $\mathcal{B}_t$ for $s\in 2C_t$.
\item $\left\vert E(\sigma_s)-E(H(\sigma_s,\frac{1}{2^{\tau-1}}\mathcal{B}_t)) - \left(  E(\sigma_s)-E(H(\sigma_s,\frac{1}{2^{\tau-1}}\mathcal{B}_t))  \right) \right\vert \leq \frac{e_{\tilde{\sigma}, \frac{\eps_0}{3^{\tau}}}(t)}{4}$
\end{itemize}
Since $A$ is compact, let $C_{t_1},\cdots,C_{t_m}$ be a Besicovitch covering of $A$ that is
$$ \forall t \in \mathbb{B}^{k-2} , 1\leq \sum_{i=1}^m \mathbf{1}_{C_{t_i}}(t) \leq \tau \hskip.1cm. $$
For $1\leq j \leq m$, we set $r_j : \mathbb{B}^{k-2} \rightarrow [0,1]$ a continuous map which satisfies $r_j = 1$ on  $C_{t_j}$ and $r_j=0$ on $\mathbb{B}^{k-2} \setminus \left(2C_{t_j}\right) \cup \bigcup_{\{i;C_{t_i}\cap C_{t_j} =\emptyset \}} C_{t_i} $.

Then, $r_j = 0$ on $C_{t_i}$ if $C_{t_i}\cap C_{t_j} =\emptyset$ so that for $t\in\mathbb{B}^{k-2}$, $r_j(t)$ is positive at most $\tau$ times. We let $\mathcal{B}_j = \mathcal{B}_{t_j}$ and then
$$\int_{\bigcup_{B\in\mathcal{B}_j} B} \left\vert \nabla \tilde{\sigma}_t \right\vert^2 \leq \frac{\eps_0}{3^{\tau-1}} \hskip.1cm.$$
Finally, for $t\in A$, there is $j(t)$ such that
$$ E(\tilde{\sigma}_t) - E\left(H\left(\tilde{\sigma}_t,\frac{r_{j(t)}(t)}{2^{\tau-1}}B_{j(t)}\right)\right) \geq \frac{e_{\tilde{\sigma}, \frac{\eps_0}{3^{\tau}}}(t_{j(t)})}{4} \geq \frac{e_{\tilde{\sigma}, \frac{\eps_0}{3^{\tau+1}}}(t)}{8}$$
which completes the proof of proposition \ref{deffamilies}.

\hfill $\diamondsuit$

Now, we can state the main theorem of the section. Notice that for the proof of this theorem, we follow the STEP 3 in the section sketch of proof, working carefully on the replacement all along the sweepout.

\begin{theo} \label{harmrepl} We fix $W>0$. There is a constant $\eps_1>0$, a constant $0<\eta\leq 1$ and a continuous map $\Psi : \mathbb{R}_+\to \mathbb{R}_+$ with $\Psi(0)=0$ such that for any sweepout $\tilde{\sigma}:\mathbb{B}^{k-2}\to \mathcal{A}$ such that 
\begin{itemize}
\item for any $t \in \mathbb{B}^{k-2}$, if $\tilde{\sigma}_t$ is a harmonic map with free boundary, then $\tilde{\sigma}_t$ is constant,
\item $\ds \max_{t} E(\sigma_t) \geq W$,
\end{itemize}
there exists a sweepout $\sigma :\mathbb{B}^{k-2}\to \mathcal{A}$ homotopic to $\tilde{\sigma}$ such that 
\begin{itemize}
\item For any $t\in\mathbb{B}^{k-2}$, $E(\sigma_t)\leq E(\tilde{\sigma}_t)$.
\item For any $t\in\mathbb{B}^{k-2}$ with $E(\tilde{\sigma}_t)\geq \frac{W}{2}$, then for any $\mathcal{B}\in\mathcal{E}$ such that the energy of $\sigma_t$ on $\bigcup_{B\in\mathcal{B}} B$ is less than $\eps_1$, then
$$ \int_{\mathbb{D}} \left\vert \nabla \left( \sigma_t - H(\sigma_t, \eta \mathcal{B}) \right) \right\vert^2  \leq \Psi ( E(\sigma_t)-E(\tilde{\sigma}_t)) \hskip.1cm.$$
\end{itemize}
\end{theo}

{\it Proof of Theorem \ref{harmrepl} :}\\

Let $\mathcal{B}_1,\cdots,\mathcal{B}_m \in \mathcal{E}$ and $r_1,\cdots,r_m$ given by proposition  \ref{deffamilies} on $\tilde{\sigma}$. We set by induction on $1\leq j\leq m$
$$ \begin{cases}
\sigma^0 = \tilde{\sigma} & \\
\sigma^j_t = H(\sigma_t^{j-1}, r_j(t) \mathcal{B}_j)  & \hbox{ for any } t\in \mathbb{B}^{k-2}
\end{cases}$$
and we set $\sigma = \sigma^m$. By proposition \ref{continuity}, it is clear that $\sigma :\mathbb{B}^{k-2}\to \mathcal{A}$ is continuous. Moreover, using again proposition \ref{continuity}, we can prove by induction that for any $1\leq j\leq m$, $\sigma^j$ is homotopic to $\sigma^0$. Indeed, the map
$$ F(\theta)_t = H(\sigma^{j-1}_t , \theta r_j(t) \mathcal{B}_j ) $$
for $\theta\in [0,1]$ and $t\in\mathbb{B}^{k-2}$ defines a homotopy between $\sigma^j$ and $\sigma^{j-1}$.\\

We fix $t \in \mathbb{B}^{k-2}$ such that $E(\tilde{\sigma}_t)\geq \frac{W}{2}$ and let $j(t)$ given by proposition  \ref{deffamilies}, such that 
\begin{equation} \label{diffengeq} E(\tilde{\sigma}_t) - E\left(H\left(\tilde{\sigma}_t,\frac{r_{j(t)}(t)}{2^{\tau-1}}B_{j(t)}\right)\right) \geq \frac{e_{\tilde{\sigma}, \frac{\eps_0}{3^{\tau+1}}}(t)}{8} \hskip.1cm. \end{equation}
Up to permutation, without changing the order in $\{ j \in\{1,\cdots,m\} ; r_j(t)>0 \}$, we assume that $r_j(t)=0$ for $j \geq \tau+1$. Therefore, $\sigma_t = H(\tilde{\sigma}_t,\mathcal{B}_1,\cdots,\mathcal{B}_{\tau})$ and $j(t)\in \{1,\cdots,\tau\}$. For $1\leq i \leq \tau$, using \eqref{exchange1} $i-1$ times, we can prove that
$$ E(\tilde{\sigma}_t) - E(H(\tilde{\sigma}_t, \frac{1}{2^{i-1}} \mathcal{B}_i)) \leq \left(1+\frac{i-1}{\kappa}\right) \left( E(\tilde{\sigma}_t) - E(H(\tilde{\sigma}_t, \mathcal{B}_1,\cdots, \mathcal{B}_i)) \right)^{\frac{1}{2}} \hskip.1cm. $$
We write this for $i = j(t)$ noticing that $j(t)\leq \tau$,
$$ E(\tilde{\sigma}_t) - E(H(\tilde{\sigma}_t, \frac{1}{2^{\tau-1}} \mathcal{B}_{j(t)})) \leq \left(1+\frac{\tau-1}{\kappa}\right) \left( E(\tilde{\sigma}_t) - E(\sigma_t) \right)^{\frac{1}{2}} \hskip.1cm. $$
With \eqref{diffengeq}, we get a constant $K$ such that
\begin{equation} \label{este}
e_{\tilde{\sigma}, \frac{\eps_0}{3^{\tau+1}}}(t) \leq K \left( E(\tilde{\sigma}_t) - E(\sigma_t) \right)^{\frac{1}{2}}
\end{equation}

Now, we set $\eps_1 = \frac{\eps_0}{3^{\tau+1}} $ and $\eta = \frac{1}{2^{2\tau-1}}$ and we let $\mathcal{B} \in\mathcal{E}$ be such that the energy of $\sigma_t$ on $\bigcup_{B\in\mathcal{B}} B$ is less than $\eps_1$. Using \eqref{exchange2} $\tau$ times, we get
$$ E(\sigma_t) - E(H(\sigma_t,\eta \mathcal{B})) \leq  E(\tilde{\sigma}_t) - E(H(\tilde{\sigma}_t, \frac{1}{2^{\tau-1}}\mathcal{B})) + \frac{\tau}{\kappa} \left( E(\tilde{\sigma}_t) - E(\sigma_t) \right)^{\frac{1}{2}} \hskip.1cm.$$
By the definition of $e_{\tilde{\sigma}, \frac{\eps_0}{3^{\tau+1}}}$ and \eqref{este}, we complete the proof of Theorem \ref{harmrepl} setting  $\Psi(r) = \left(K + \frac{\tau}{\kappa}\right) r^{\frac{1}{2}}$.

\hfill $\diamondsuit$

\section{Energy identity}

In this section, we aim at proving theorem \ref{palaissmaleassump}. The min-max sequences satisfy the assumptions \eqref{epsregesttildeuj} and \eqref{quasiconfuj} of the theorem so that we get the energy identity. The bubble tree construction is now classical as soon as we are able to prove a fundamental no-neck-energy lemma. We do not detail this construction since it is standard, written in \cite{LP} for instance.

Let's focus on the no-neck-energy lemma. It was already proved in the interior case in proposition \cite{CM2008}. Here we focus mostly in the boundary case, where the main improvement is in proposition \ref{noneckangularenergy} where thanks to our previous work \cite{LP}, see also \cite{LR}, we are able to prove that the angular part of the gradient is small as soon as the energy is small enough.

\begin{prop} \label{noneckenergy}
There is a constant $\eps_2>0$ such that for any $0<\eta\leq 1$, for any sequence $\lambda_n \to 0$, for any sequence of maps $\{u_n\}$ in $W^{1,2}(\mathbb{D}_+ \setminus \mathbb{D}_{\lambda_n}, N)$ satisfying $u_n( I \setminus \mathbb{D}_{\lambda_n} ) \in M$ such that

\begin{equation*}
\begin{split}
&\int_{\mathbb{D}_+ \setminus \mathbb{D}_{\lambda_n}} \left\vert \nabla u_n \right\vert^2 \leq \eps_2,\\
& \sup_{\underset{\forall B\in \mathcal{B}, B\in \mathbb{D}_+ \setminus \mathbb{D}_{\lambda_n}}{\mathcal{B} \in \mathcal{E} \hbox{ \small{and} }} } \left\{ \int_{\mathbb{D}_+ \setminus \mathbb{D}_{\lambda_n}} \left\vert \nabla(u_n - H(u_n,\eta \mathcal{B}))\right\vert^2  \right\}  = o(1)
\end{split}
\end{equation*}
and 
$$\frac{1}{2}\int_{\mathbb{D}_+ \setminus \mathbb{D}_{\lambda_n}} \left\vert \nabla u_n \right\vert^2  - \int_{\mathbb{D}_+ \setminus \mathbb{D}_{\lambda_n}} \sqrt{\left\vert \nabla_{\theta} u_n \right\vert^2\left\vert \nabla_r u_n \right\vert^2 - \left\langle \nabla_{\theta} u_n, \nabla_r u_n \right\rangle^2 } = o(1) \hskip.1cm.$$

Then,
$$ \lim_{R\to +\infty} \limsup_{n\to +\infty} \int_{\mathbb{D}_+ \cap \left(\mathbb{D}_{\frac{1}{R}} \setminus \mathbb{D}_{R \lambda_n}\right)} \left\vert \nabla u_n \right\vert^2 = 0 \hskip.1cm. $$

\end{prop}

In fact, thanks to the first assumption and the second assumption we aim at proving that
\begin{equation}\label{noneckangularenergy} \lim_{R\to +\infty} \limsup_{n\to +\infty} \int_{\mathbb{D}_+ \cap \left(\mathbb{D}_{\frac{1}{R}} \setminus \mathbb{D}_{R \lambda_n}\right)} \left\vert \nabla_{\theta} u_n \right\vert^2 = 0 \hskip.1cm. \end{equation}
Indeed, we easily deduce from the third assumption that for any $R>0$, 
\begin{equation} \label{angulareqradial} \int_{\mathbb{D}_+ \cap \left(\mathbb{D}_{\frac{1}{R}} \setminus \mathbb{D}_{R \lambda_n}\right)} \left\vert \nabla_{r} u_n \right\vert^2 = \int_{\mathbb{D}_+ \cap \left(\mathbb{D}_{\frac{1}{R}} \setminus \mathbb{D}_{R \lambda_n}\right)} \left\vert \nabla_{\theta} u_n \right\vert^2 + o(1) \hbox{ as } n\to +\infty  \end{equation}
so that once \eqref{noneckangularenergy} is proved, we get proposition \ref{noneckenergy}. Notice that if $\{u_n\}$ is already a sequence of harmonic maps defined on $\mathbb{D}_+$ with free boundary in $M$ on $I$, the third assumption can be deleted since \eqref{angulareqradial} is already true by a Poho\v{z}aev identity, see \cite{LP} for the case $M=\mathbb{S}^n$ which extend straightforward to general target. However, the third assumption is in some sense stronger because in this case, we do not need $u_n$ to be defined on $\overline{\mathbb{D}_+}$.

In order to prove \eqref{noneckangularenergy}, we will first prove a weaker property than \eqref{noneckangularenergy}, for free boundary harmonic maps.

\begin{prop} \label{angularenergyharm}
For any $\delta>0$, there exists $\eps_3>0$ and $\alpha>1$ such that for any harmonic map $u : \mathbb{D}_+ \setminus \mathbb{D}_{\frac{1}{\alpha^3}} \to N$ with free boundary in $M$ on $I \setminus \mathbb{D}_{\frac{1}{\alpha^3}}$ satisfying 
$$\int_{\mathbb{D}_+ \setminus \mathbb{D}_{\frac{1}{\alpha^3}}} \left\vert \nabla u \right\vert^2 \leq \eps_3 \hskip.1cm,$$
we have 
$$ \int_{\mathbb{D}_+ \cap \left( \mathbb{D}_{\frac{1}{\alpha}} \setminus \mathbb{D}_{\frac{1}{\alpha^2}}\right)} \left\vert \nabla_\theta u \right\vert^2 < \delta \int_{\mathbb{D}_+ \setminus \mathbb{D}_{\frac{1}{\alpha^3}}} \left\vert \nabla u \right\vert^2 \hskip.1cm.$$
\end{prop}

{\it Proof of proposition \ref{angularenergyharm} :}\\

Let $0<\eps_3<\eps_0$, with $\eps_0$ as in Proposition \ref{eereg} and $\eps_3$ will be fixed later. As in the proof of Proposition \ref{eereg}, let us consider $\tilde{u}$ the extension of $u$ on $\mathbb{D} \setminus \mathbb{D}_{\frac{1}{\alpha^3}}$. A priori, we have to reduce the ball, but since the final result is on a reduced ball, we don't do it not to make the notation heavier. As in Claim 3.1 of \cite{LP}\footnote{In fact the theorem deals with free boundary harmonic maps into the ball, but after extension the equation is exactly the same as here, hence it can be applied verbatim.} or the seminal paper \cite{LR}, we can choose $\alpha_0>1$ such that, 

$$\Vert \nabla_\theta \tilde{u}\Vert_{L^{2,1}\left(\mathbb{D}_{\frac{1}{\alpha_0}} \setminus \mathbb{D}_{\frac{1}{\alpha_0^2}}\right)} \leq C \int_{ \D\setminus \mathbb{D}_{\frac{1}{\alpha_0^3}}} \left\vert \nabla \tilde{u} \right\vert^2 \leq C \int_{ \D_+\setminus \mathbb{D}_{\frac{1}{\alpha_0^3}}} \left\vert \nabla u \right\vert^2, $$
the last inequality comes from the fact that we uniformly control the $L^2$-norm of the extension by the one of the initial map, see Proposition \ref{eereg}. Then it suffices to prove that $L^{2,\infty}$-norm of $\nabla \tilde{u}$ can be made as small as we want up to take $\eps_3<\eps_0$ small enough and $\alpha >\alpha_0$ big enough.  For this it suffices to prove that for any $\eta>0$, there exists $\eps_3<\eps_0$ and $\alpha>\alpha_0$ such that for all $r\in \left(\frac{4}{\alpha^2},\frac{1}{4\alpha}\right)$, we have
$$
\int_{\D_{2r }\setminus \mathbb{D}_{\frac{r}{2}}} \left\vert \nabla \tilde{u} \right\vert^2 \leq \eta \int_{ \D_+\setminus \mathbb{D}_{\frac{1}{\alpha^3}}} \left\vert \nabla u \right\vert^2.
$$
Indeed, by $\eps$-regularity we have for all $x\in \D_\frac{1}{2\alpha}\setminus \D_\frac{2}{\alpha^2}$ that
$$\vert \nabla \tilde{u}(x) \vert \leq \frac{C}{\vert x\vert}  \sqrt{\int_{ \D_{2\vert x\vert}\setminus \mathbb{D}_{\frac{\vert x\vert}{2}}} \left\vert \nabla \tilde{u} \right\vert^2}. $$
Then using the fact that $\frac{1}{\vert x\vert}$ is into $L^{2,\infty}(\R^2)$, we will get the result.\\

Now, let us assume by contradiction, that there exists a sequence $u_n$ satisfying the assumptions of the theorem and a sequence $r_n \rightarrow 0 $ such that  

\begin{equation}
\label{ABL}
\int_{ \D_{2r_n}\setminus \mathbb{D}_{r_n/2}} \left\vert \nabla \tilde{u}_n \right\vert^2 >\eta \int_{ \D_+\setminus \mathbb{D}_{\frac{1}{\alpha_n^3}}} \left\vert \nabla u_n \right\vert^2 \hskip.1cm.
\end{equation}
and
\begin{equation}
\label{ABL}
\int_{ \D_+\setminus \mathbb{D}_{\frac{1}{\alpha_n^3}}} \left\vert \nabla u_n \right\vert^2=\eps_n\rightarrow 0 \hbox{ as } \alpha_n\rightarrow +\infty\hskip.1cm.
\end{equation}

Then, setting $v_n=\frac{\tilde{u}_n\left(\, .\, r_n\right)-\fint_{\partial\D}\tilde{u}_n\left(\, .\, r_n\right)}{\sqrt{\eps_n}}$ and using once again the $\eps$-regularity, we easily see that $v_n$ converges in $W^{1,2}_{loc}(\R^2\setminus\{0\})$ to a non trivial harmonic map which takes values in some tangent space identifies with some $\R^n$, since obtained from a free boundary harmonic map by reflection. Still by $\eps$-regularity we have that there exists $C>0$ such that 

$$\vert \nabla v(x) \vert \leq \frac{C}{\vert x\vert},$$

Thanks to the classical B\^{o}cher theorem \cite{SBR}, we necessarily get that 
$$v(x)=a\ln(\vert x\vert)+b.$$
Finally, remembering that $\nabla v$ must be bounded in $L^{2}(\R^2\setminus\{0\})$ by conformal invariance, $v$ must be constant, which is a contradiction with \eqref{ABL}. This achieves the proof of the proposition.

\hfill $\diamondsuit$\\

Now, using proposition \ref{angularenergyharm}, we prove a similar result for a more general class of maps than harmonic maps with free boundary, the conclusion is still weaker than \eqref{noneckangularenergy}.

\begin{prop} \label{angularenergyalmostharm}
For any $\delta>0$, there exists $\eps_4>0$ and $\alpha>1$ such that for any $\eta>0$, there exists $\mu>0$ such that for any map $u \in W^{1,2}( \mathbb{D}_+ \setminus \mathbb{D}_{\frac{1}{\alpha^5}} , N)$ with $u(I\setminus \mathbb{D}_{\frac{1}{\alpha^5}})\subset M$ satisfying
$$\int_{\mathbb{D}_+ \setminus \mathbb{D}_{\frac{1}{\alpha^5}}} \left\vert \nabla u \right\vert^2 \leq \eps_4 $$
and
$$ \sup_{\underset{\forall B\in \mathcal{B}, B\subset \mathbb{D}_+ \setminus \mathbb{D}_{\frac{1}{\alpha^5}}}{\mathcal{B} \in \mathcal{E} \hbox{ \small{and} }} } \left\{ \int_{\mathbb{D}_+ \setminus \mathbb{D}_{\frac{1}{\alpha^5}}} \left\vert \nabla(u - H(u,\eta \mathcal{B}))\right\vert^2  \right\} \leq \mu \int_{\mathbb{D}_+ \setminus \mathbb{D}_{\frac{1}{\alpha^5}}} \left\vert \nabla u \right\vert^2 \hskip.1cm,$$
we have
$$ \int_{\mathbb{D}_+ \cap \left( \mathbb{D}_{\frac{1}{\alpha^2}} \setminus \mathbb{D}_{\frac{1}{\alpha^3}}\right)} \left\vert \nabla_\theta u \right\vert^2 \leq \delta \int_{\mathbb{D}_+ \setminus \mathbb{D}_{\frac{1}{\alpha^5}}} \left\vert \nabla u \right\vert^2 \hskip.1cm.$$
\end{prop}

{\it Proof of proposition \ref{angularenergyalmostharm} :}\\

Let $\eps_3>0$ and $\alpha>0$ given by proposition \ref{angularenergyharm}. By contradiction, we assume that there is a sequence of maps $u_n \in W^{1,2}( \mathbb{D}_+ \setminus \mathbb{D}_{\frac{1}{\alpha^5}} , N)$ with $u_n(I\setminus \mathbb{D}_{\frac{1}{\alpha^5}})\subset M$ satisfying
\begin{equation} \label{smallen} \int_{\mathbb{D}_+ \setminus \mathbb{D}_{\frac{1}{\alpha^5}}} \left\vert \nabla u_n \right\vert^2 \leq \eps_3 \hskip.1cm,\end{equation}
\begin{equation} \label{almostharm} \sup_{\underset{\forall B\in \mathcal{B}, B\in \mathbb{D}_+ \setminus \mathbb{D}_{\frac{1}{\alpha^5}}}{\mathcal{B} \in \mathcal{E} \hbox{ \small{and} }} } \left\{ \int_{\mathbb{D}_+ \setminus \mathbb{D}_{\frac{1}{\alpha^5}}} \left\vert \nabla(u_n - H(u_n,\eta \mathcal{B}))\right\vert^2  \right\} \leq \frac{1}{n} \int_{\mathbb{D}_+ \setminus \mathbb{D}_{\frac{1}{\alpha^5}}} \left\vert \nabla u_n \right\vert^2 \hskip.1cm, \end{equation}
and
\begin{equation} \label{absurdconclusion} \int_{\mathbb{D}_+ \cap \left( \mathbb{D}_{\frac{1}{\alpha^2}} \setminus \mathbb{D}_{\frac{1}{\alpha^3}}\right)} \left\vert \nabla_\theta u_n \right\vert^2 > \delta \int_{\mathbb{D}_+ \setminus \mathbb{D}_{\frac{1}{\alpha^5}}} \left\vert \nabla u_n \right\vert^2 \hskip.1cm. \end{equation}

\hspace{4mm}

\textit{CASE 1 :} We assume that $\int_{\mathbb{D}_+ \setminus \mathbb{D}_{\frac{1}{\alpha^5}}} \left\vert \nabla u_n \right\vert^2$ is bounded from below. 

Then, with \eqref{smallen} and \eqref{almostharm}, it is clear that $u_n$ converges in $W^{1,2}\left( \mathbb{D}_+ \cap \left(\mathbb{D}_{\frac{1}{\alpha}}\setminus \mathbb{D}_{\frac{1}{\alpha^4}}\right)  , N\right)$ to some a harmonic map $u : \mathbb{D}_+ \cap \left(\mathbb{D}_{\frac{1}{\alpha}}\setminus \mathbb{D}_{\frac{1}{\alpha^4}}\right) \to N$ with free boundary in $M$ on $I \cap \left(\mathbb{D}_{\frac{1}{\alpha}}\setminus \mathbb{D}_{\frac{1}{\alpha^4}}\right)$. Then passing to the limit in \eqref{absurdconclusion} gives a contradiction, by proposition \ref{angularenergyharm}.

\hspace{4mm}

\textit{CASE 2 :} Up to a subsequence, we assume that $\int_{\mathbb{D}_+ \setminus \mathbb{D}_{\frac{1}{\alpha^5}}} \left\vert \nabla u_n \right\vert^2 \to 0$. Then, $u_n \to p \in M$ in $W^{1,2}$ as $n\to +\infty$. Then, we are going to 
blow-up $u_n$ arround $p$ in order to produce a harmonic map into $T_pN$ with free boundary into $T_p M$.
 
Take a family of balls $\{B_i\}_{1\leq i\leq m}$ centered in $x_i \in \mathbb{D}_{\frac{1}{\alpha}}^+\setminus \mathbb{D}_{\frac{1}{\alpha^4}}$ and half balls centered in $x_i \in I \cap \left(\mathbb{D}_{\frac{1}{\alpha}}\setminus \mathbb{D}_{\frac{1}{\alpha^4}}\right)$ included in  $\mathbb{D} \setminus \mathbb{D}_{\frac{1}{\alpha^5}}$ such that 
$$\mathbb{D}_{+} \cap \left( \mathbb{D}_{\frac{1}{\alpha}}\setminus \mathbb{D}_{\frac{1}{\alpha^4}}\right) \subset \bigcup_{i=1}^k \frac{\eta}{2}B_i \hskip.1cm.$$
We set for $1\leq i \leq m$
$$ v_n^i = H(u_n,\eta B_i) \hbox{ and } \tilde{v}_n^i = \left(v_n^i - \overline{v_n^i} \right)\left( \int_{\mathbb{D}_+ \setminus \mathbb{D}_{\frac{1}{\alpha^5}}} \left\vert \nabla u_n \right\vert^2\right)^{-\frac{1}{2}}  \hbox{ and }$$
$$ \tilde{u}_n = \left( u_n - \overline{u_n} \right) \left(\int_{\mathbb{D}_+ \setminus \mathbb{D}_{\frac{1}{\alpha^5}}} \left\vert \nabla u_n \right\vert^2\right)^{-\frac{1}{2}} $$
where $\overline{u_n}$ is the mean of $u_n$ on $\mathbb{D}_+ \setminus \mathbb{D}_{\frac{1}{\alpha^5}}$. For $1\leq i \leq m$, it is clear that up to subsequences,
$$ \tilde{v}_n^i \to \tilde{v}^i  \in \mathcal{C}^2\left(\frac{\eta}{2}B_i\right) \hbox{ as } n \to +\infty$$
where $\tilde{v}^i$ is harmonic into $T_pN$ and $\partial_{\nu}  \tilde{v}_i(x) \in (T_{p}M)^{\perp} $ for any $x \in I \cap \left(\frac{\eta}{2}B_i\right)$. Now, by \eqref{almostharm},
$$ \int_{\eta B_i} \left\vert\nabla (\tilde{u}_n - \tilde{v}_n^i )\right\vert^2 \leq \frac{1}{n}$$
so that $\nabla \tilde{u}_n \to \nabla \tilde{v}^i$ in $L^2\left(\frac{\eta}{2}B_i\right)$. Moreover, by a Poincar\'e inequality, $\tilde{u}_n$ is also bounded in $L^2$ so that it converges to some map $\tilde{u}$ in $W^{1,2}\left(\mathbb{D}_{+} \cap \left( \mathbb{D}_{\frac{1}{\alpha}}\setminus \mathbb{D}_{\frac{1}{\alpha^4}}\right)\right)$. On $\frac{\eta}{2}B_i$, this map is equal to $ \tilde{v}^i$ up to constant so that $\tilde{u}$ is also harmonic in $T_pN$ with free boundary in $T_pM$ on $I\cap \left( \mathbb{D}_{\frac{1}{\alpha}}\setminus \mathbb{D}_{\frac{1}{\alpha^4}}\right)$.  Passing to the limit in \eqref{absurdconclusion}, we get 
\begin{equation}  \int_{\mathbb{D}_+ \cap \left( \mathbb{D}_{\frac{1}{\alpha^2}} \setminus \mathbb{D}_{\frac{1}{\alpha^3}}\right)} \left\vert \nabla_{\theta} \tilde{u}_n \right\vert^2 \geq \delta \hskip.1cm, \end{equation}
which contradicts proposition \ref{angularenergyharm}, since this last one is true for any $\eps_3$ for harmonic maps into $T_p N \cong \R^n$ with free boundary into $T_pM\cong \R^m$.

\hfill $\diamondsuit$

Now, we are able to complete the proof.\\

{\it Proof of proposition \ref{noneckenergy} :}\\

Thanks to the remark after Proposition \ref{noneckenergy}, we just have to prove \eqref{noneckangularenergy}. We fix $\delta >0$, let $\alpha>1$, $\eps_2=\eps_4$ and $\mu >0$ be given by proposition \ref{angularenergyalmostharm} and let $R > \alpha^3$. We aim at proving that there is some constant $C$ such that 
\begin{equation}\label{purpose} \limsup_{n\to +\infty} \int_{\mathbb{D}_+ \cap \left(\mathbb{D}_{\frac{1}{R}} \setminus \mathbb{D}_{R \lambda_n}\right)} \left\vert \nabla_{\theta} u_n \right\vert^2 \leq C\delta \hskip.1cm. \end{equation}
We let $N_n$ be the smallest integer greater than $\frac{-2\ln(R)-\ln(\lambda_n)}{\ln{\alpha}}$, we set
$$ \mathbb{A}_k = \mathbb{D}_+ \cap \left( \mathbb{D}_{\frac{1}{R\alpha^k}} \setminus \mathbb{D}_{\frac{1}{R\alpha^{k+1}}} \right) \hbox{ and }  \tilde{\mathbb{A}}_k = \mathbb{D}_+ \cap \left( \mathbb{D}_{\frac{1}{R\alpha^{k-2}}} \setminus \mathbb{D}_{\frac{1}{R\alpha^{k+3}}} \right) $$
 for $0\leq k \leq N_n -1$, and
$$ E_n = \{ k\in[0, N_n -1] ;  \sup_{\underset{\forall B\in \mathcal{B}, B\subset \tilde{\mathbb{A}}_k }{\mathcal{B} \in \mathcal{E} \hbox{ \small{and} }} } \left\{ \int_{\tilde{\mathbb{A}}_k} \left\vert \nabla(u_n - H(u_n,\eta \mathcal{B}))\right\vert^2  \right\}  \leq \mu \int_{\tilde{\mathbb{A}}_k} \left\vert \nabla u_n \right\vert^2  \} \hskip.1cm,$$
$$ F_n = \{ k\in[0, N_n -1] ;  \sup_{\underset{\forall B\in \mathcal{B}, B\subset \tilde{\mathbb{A}}_k }{\mathcal{B} \in \mathcal{E} \hbox{ \small{and} }} } \left\{ \int_{\tilde{\mathbb{A}}_k} \left\vert \nabla(u_n - H(u_n,\eta \mathcal{B}))\right\vert^2  \right\}  > \mu \int_{\tilde{\mathbb{A}}_k} \left\vert \nabla u_n \right\vert^2  \}$$
so that 
\begin{equation} \label{sumEF} \int_{\mathbb{D}_+ \cap \left(\mathbb{D}_{\frac{1}{R}} \setminus \mathbb{D}_{R \lambda_n}\right)} \left\vert \nabla_{\theta} u_n \right\vert^2 = \sum_{k\in E_n} \int_{\mathbb{A}_k} \left\vert \nabla_{\theta} u_n \right\vert^2 + \sum_{k\in F_n} \int_{\mathbb{A}_k} \left\vert \nabla_{\theta} u_n \right\vert^2 \hskip.1cm. \end{equation}
By proposition \ref{angularenergyalmostharm},
\begin{equation}\label{sumE} \sum_{k\in E_n} \int_{\mathbb{A}_k} \left\vert \nabla_{\theta} u_n \right\vert^2 \leq \sum_{k\in E_n} \delta \int_{\tilde{\mathbb{A}}_k}  \left\vert \nabla u_n \right\vert^2 \leq \delta \sum_{k=0}^{N_n -1} \sum_{i=k-2}^{k+2}  \int_{\mathbb{A}_i}  \left\vert \nabla u_n \right\vert^2 \leq 5 \eps_2 \delta \hskip.1cm. \end{equation}
Now let $k\in F_n$ and $\mathcal{B}_n^k \in \mathcal{E}$ with $\forall B\in \mathcal{B}_n^k$, $B \subset \tilde{\mathbb{A}}_k$ such that 
$$ \left\{ \int_{\tilde{\mathbb{A}}_k} \left\vert \nabla(u_n - H(u_n,\eta \mathcal{B}_n^k))\right\vert^2  \right\}  > \mu \int_{\tilde{\mathbb{A}}_k} \left\vert \nabla u_n \right\vert^2 \hskip.1cm.$$
We set for $1\leq j \leq 5$
$$ \tilde{\mathcal{B}}_n^j = \bigcup_{ \underset{k = j \mod 5}{k\in F_n} } \mathcal{B}_n^k \hskip.1cm.$$
Then, $\tilde{\mathcal{B}}_n^j \in \mathcal{E}$ and there is $1 \leq j \leq 5$ such that
\begin{eqnarray*}
 \int_{\mathbb{D}_+\setminus \mathbb{D}_{\lambda_n}} \left\vert \nabla(u_n - H(u_n,\eta \tilde{\mathcal{B}}_n^j))\right\vert^2 & \geq & \frac{1}{5} \sum_{i=1}^{5} \int_{\mathbb{D}_+\setminus \mathbb{D}_{\lambda_n}} \left\vert \nabla(u_n - H(u_n,\eta \tilde{\mathcal{B}}_n^i))\right\vert^2 \\
 & \geq & \frac{1}{5} \sum_{i=1}^{5}  \sum_{ \underset{k = i \mod 5}{k\in F_n} }  \int_{\mathbb{D}_+\setminus \mathbb{D}_{\lambda_n}} \left\vert \nabla(u_n - H(u_n,\eta \mathcal{B}_n^k))\right\vert^2 \\
 & \geq & \frac{\mu}{5} \sum_{k\in F_n} \int_{\tilde{A}_k} \left\vert \nabla u_n \right\vert^2 \\
 & \geq & \frac{\mu}{5} \sum_{k\in F_n} \int_{A_k} \left\vert \nabla u_n \right\vert^2 \, .
\end{eqnarray*} 
Then using the second assumption in proposition \ref{noneckenergy}, we get
$$  \sum_{k\in F_n} \int_{A_k} \left\vert \nabla u_n \right\vert^2 = o(1) \hbox{ as } n\to +\infty \hskip.1cm. $$
Combining this with \eqref{sumEF} and \eqref{sumE}, we get \eqref{purpose} with $C= 5 \eps_2$. This completes the proof of \eqref{noneckangularenergy}, and as already said Proposition \ref{noneckenergy} thanks to the remark after it.

\hfill $\diamondsuit$

\section{Proof of Theorem \ref{main}}

In this section, we aim at gathering all the previous results in order to prove the main theorem. Before applying the replacement procedure in theorem \ref{harmrepl}, we need to build a minimizing sequence of sweepouts satisfying that any free boundary harmonic slice has to be constant. The proof is in the same spirit as Theorem 2.1 in \cite{CM2008} 

\begin{prop} \label{harmonicthenconstant}
Let $\sigma : \mathbb{B}^{k-2} \to \mathcal{A}$ be a sweepout and $\eps >0$. Then, there is a sweepout $\tilde{\sigma} : \mathbb{B}^{k-2} \to \mathcal{A}$ homotopic to $\sigma$ such that
\begin{itemize}
\item $\forall t \in \mathbb{B}^{k-2}, \Vert \tilde{\sigma}_t - \sigma_t \Vert_{W^{1,2}} \leq \eps$
\item $\forall t \in \mathbb{B}^{k-2}$ if $\tilde{\sigma}_t: \mathbb{D}\to N$ is a harmonic map with free boundary in $M$, then $\tilde{\sigma}_t$ is a constant map.
\end{itemize}
\end{prop}

{\it Proof of Proposition \ref{harmonicthenconstant} :}\\

We first set for $0<\eta<\frac{1}{2}$ and $x\in \mathbb{D}_{1-\eta}$
$$ \sigma_t^{\eta} = \int_{\mathbb{R}^2} \frac{\phi\left(\frac{y}{\eta}\right)}{\eta^2} \sigma_t(x-y)dy $$
a regularization of $\sigma_t^{\eta}$ where $\phi \geq 0$ is a smooth function with compact support in $\mathbb{D}$. By classical results on convolutions, we have that $\sigma_t^{\eta} \in \mathcal{C}^{\infty}(\mathbb{D}_{\frac{1}{2}})$, that $t\mapsto \sigma_t^{\eta}$ is continuous from $\mathbb{B}^{k-2}$ to $\mathcal{C}^2(\mathbb{D}_{\frac{1}{2}})$ and
\begin{equation}
\label{approxconvol} \lim_{\eta\to 0} \sup_{t\in \mathcal{B}^{k-2}} \Vert \sigma_t^{\eta} - \sigma_t \Vert_{W^{1,2}\cap\mathcal{C}^0\left(\mathbb{D}_{\frac{1}{2}}\right)} = 0 \hskip.1cm.
\end{equation}
We can define for $\eta$ small enough
$$ \tilde{\sigma}_t^{\eta}(x) = \pi_N(\theta(x)\sigma_t^{\eta}(x) + (1-\theta(x))\sigma_t^{\eta}(x)) $$
where $\theta \in \mathcal{C}_c^{\infty}(\mathbb{D}_{\frac{1}{2}})$, $0 \leq \theta\leq 1$, $\theta=1$ in $\mathbb{D}_{\frac{1}{4}}$. Let $\eta>0$ be such that
\begin{equation}
\label{approxconvol2} \sup_{t\in \mathbb{B}^{k-2}} \Vert \tilde{\sigma}_t^{\eta} - \sigma_t \Vert_{W^{1,2}\cap\mathcal{C}^0\left(\mathbb{D}\right)} \leq \frac{\eps}{2} \hskip.1cm.
\end{equation}
Now, we define for $0<\rho<\frac{1}{4}$ a map $\Phi_{\rho} : \mathbb{R}^2 \to \mathbb{R}^2$ by
$$ \Phi_{\rho}(r\cos \theta,r\sin\theta) = \begin{cases} 2(r\cos \theta,r\sin\theta) \hbox{ if } r\leq \frac{\rho}{2} \\
 (\rho \cos \theta,\rho \sin\theta) \hbox{ if } \frac{\rho}{2}\leq r\leq \rho \\
(r\cos \theta,r\sin\theta) \hbox{ if } r\geq \rho 
 \end{cases} \hskip.1cm.$$
and we set
$$ \tilde{\sigma}_t = \tilde{\sigma}_t^{\eta} \circ \Phi_{\rho} $$
and let $\rho >0$ small enough such that
\begin{equation}
\label{approxconvol3} \sup_{t\in \mathbb{B}^{k-2}} \Vert \tilde{\sigma}_t^{\eta} - \tilde{\sigma}_t \Vert_{W^{1,2}\cap\mathcal{C}^0\left(\mathbb{D}\right)} \leq \frac{\eps}{2} \hskip.1cm.
\end{equation}
Then, $\tilde{\sigma}_t \in  W^{1,2}\cap \mathcal{C}^0(\mathbb{\overline{D}}, N)$ and it is clear that $\tilde{\sigma}$ is homotopic to $\sigma$ : we first define an homotopy from $\tilde{\sigma}$ to $\tilde{\sigma}^{\eta}$ by shrinking $\rho$ to $0$ and then an homotopy from $\tilde{\sigma}^{\eta}$ to $\sigma$ by shrinking $\eta$ to $0$ (with \eqref{approxconvol}). Moreover if for some $t\in \mathbb{B}^{k-2}$, $\tilde{\sigma}_t : \mathbb{D} \to N$ is harmonic with free boundary, then, it is conformal\footnote{ It can be easily seen considering teh Hopf differential of the immersion and noting it must be real on the boundary, see \cite{DHS}.}, and since $\partial_r \tilde{\sigma}_t = D\tilde{\sigma}_t^{\eta}(\Phi_{\rho}).\partial_r \Phi_{\rho} = 0 $ in $\mathbb{D}_{\rho}\setminus \mathbb{D}_{\frac{\rho}{2}}$, we get that $\tilde{\sigma}_t$ is constant in $\mathbb{D}_{\rho}\setminus \mathbb{D}_{\frac{\rho}{2}}$. Since $\tilde{\sigma}_t$ is harmonic, $ \tilde{\sigma}_t $ is constant everywhere, by classical unique continuation results for second order elliptic PDEs. Thanks to \eqref{approxconvol2} and \eqref{approxconvol3}, we have the expected approximation.

\hfill $\diamondsuit$

{\it Proof of Theorem \ref{main} :}\\

Let $\{s^n\}$ be a minimizing sequence of homotopic sweepouts for $W$, that is
$$ \max_{t\in \mathbb{B}^{k-2}} E(s_t^n) \to W \hbox{ as } n\to +\infty \hskip.1cm. $$
Indeed recalling STEP 0 of the sketch of the proof in section 2, the Min-Max for the Energy or the Area achieves the same value.\\ 
Applying Proposition \ref{harmonicthenconstant}, there is a sequence $\{\tilde{s}^n\}$ such that if for some $t\in \mathbb{B}^{k-2}$, $\tilde{s}_t^n$ is harmonic with free boundary, then $\tilde{s}_t^n$ is constant and that
$$ \sup_{t\in \mathbb{B}^{k-2}} \Vert \tilde{s}_t^n - s_t^n \Vert_{W^{1,2}} \to 0 \hbox{ as } n\to +\infty \hskip.1cm.$$
It is clear that $\{\tilde{s}^n\}$ is also a minimizing sequence for $W$. Let $\sigma^n$ be the sweepout given by the procedure of theorem \ref{harmrepl} starting from $\tilde{s}^n$. Using theorem \ref{harmrepl}, we have
\begin{equation} \label{reduceenergy} \forall t \in \mathbb{B}^{k-2}, E(\sigma_t^n) \leq E(\tilde{s}_t^n) \hskip.1cm.  \end{equation}
In particular, $\{\tilde{\sigma}^n\}$ is also a minimizing sequence for $W$. Now let $t_n$ be a sequence of parameters such that $Area(\sigma_{t_n}^n)\to W$ as $n \to +\infty$. Then
\begin{equation}
\label{qconf}
 W+o(1) = Area(\sigma^n_{t_n}) \leq E(\sigma^n_{t_n}) \leq E(\tilde{s}^n_{t_n}) \leq W +o(1) \hskip.1cm.
 \end{equation}
Then, $E(\sigma^n_{t_n}) - E(\tilde{s}^n_{t_n}) = o(1)$ as $n\to +\infty$ and by theorem \ref{harmrepl}, there exist $\eps_1>0$ and  $0< \eta\leq 1$ such that for any $\mathcal{B}\in \mathcal{E}$ such that the energy on $\bigcup_{B\in \mathcal{B}} B$ is less than $\eps_1$, we have
\begin{equation}
\label{ang1}
 \int_{\mathbb{D}} \vert \nabla \sigma_t^n - H(\sigma_t^n,\eta \mathcal{B}) \vert^2 \to 0 \hbox{ as } n\to +\infty \hskip.1cm.
 \end{equation}
We also have, thanks to \eqref{qconf}, that 
\begin{equation}
\label{ang2}
\lim_{n\to +\infty} E(\sigma^n_{t_n}) - Area(\sigma^n_{t_n}) = o(1) \hskip.1cm.
\end{equation}
Thanks to \eqref{ang1} and \eqref{ang2} and Proposition \ref{noneckenergy} we have the no-neck-energy for the sequence $\sigma^n_{t_n}$. Now we can conclude by the classical bubble tree decomposition to get the $W^{1,2}$-bubble convergence, following for instance verbatim section 3 of \cite{LP}, Step 1 and Step 2 being consequences of the fact that there is a free boundary harmonic map $W^{1,2}$-close which satisfies the $\eps$-regularity. Finally, in Step 3, Claim 3.1 has to be replace by the no-neck-energy we have just proved.\hfill $\diamondsuit$

\appendix

\section{ $\eps$-regularity for free boundary harmonic maps}
In this section we proved a generalisation of our preceding $\eps$-regularity result (see Claim 2.3 of \cite{LP}) to a general target manifold. In order to do so, we will prove some  extension of  the theorem 3.4 of \cite{JLM}.

\begin{prop} 
\label{epsreg} There is $\eps_0>0$ and a constant $C_k$ such that if a weakly harmonic map with free boundary $u\in W^{1,2}(\mathbb{D}^+,N)$  satisfies
$$ \int_{\mathbb{D}^+} \left| \nabla u\right|^2 \leq \eps_0 \hskip.1cm,$$
then $\tilde{u} \in \mathcal{C}^{\infty}(\mathbb{D}_{\frac{1}{2}}^+)$ and for any $k\geq 1$,
\begin{equation}
\label{eereg}
\left\| \nabla u \right\|_{\mathcal{C}^{k}(\mathbb{D}_{\frac{1}{2}}^+)} \leq C_k \left\| \nabla u \right\|_{L^{2}(\mathbb{D}^+)}\hskip.1cm.
\end{equation}
\end{prop}

{\it Proof :}\\
We only prove some $W^{1,p}$-estimate for $p>2$, since it is classical, see for instance \S 4 of \cite{GJ}, that once we get some $C^{0,\alpha}$ estimate we can bootstrap it, since the equation becomes sub-critical.\\

We are under the assumption of the Proposition 3.3 of \cite{JLM}, then there exists an extension of $u$,  $\tilde{u}\in W^{1,2}( \mathbb{D}_{\frac{3}{4}})$ a weak solution of 
$$div (Q \nabla \tilde{u}) =\Omega Q\nabla \tilde{u} ,$$
where  $Q\in W^{1,2} (\mathbb{D}_{\frac{3}{4}}, GL(m)) \cap L^\infty(\mathbb{D}_{\frac{3}{4}}, GL(m)) $ and $\Omega \in L^{2} (\mathbb{D}_{\frac{3}{4}}, o(m))$, moreover there exists $C> 0$ such that 
$$\vert \Omega \vert \leq C \vert \nabla \tilde{u}\vert \hbox{ almost everywhere .}$$ 

In order to prove the $\eps$-regularity, we follow the strategy of Rivi\`{e}re \cite{Riv08}, see also \cite{Schi}. Since $\Omega$ is anti symmetric, there exists $P \in W^{1,2}(\mathbb{D}_{\frac{3}{4}}, SO(m))$, such that 

$$div(~^{t}P \nabla P-~^{t}P \Omega P) =0,$$
and
$$\Vert \nabla P \Vert_2 \leq 2 \Vert \Omega \Vert_2 .$$

Hence we have ,

$$div ( \tilde{Q} \nabla \tilde{u} ) = \nabla^\bot B (\tilde{Q} \nabla \tilde{u})$$
with $\tilde{Q}=~^{t} P Q$ and $\nabla^\bot B=- ~^{t}P \nabla P+ ~^{t}P \Omega P$. Then we would like to rewrite the system like a Jacobian on the right hand-side. Let $A\in W^{1,2}(\mathbb{D}_{\frac{3}{4}}, Gl(m))$ and $C\in W^{1,2}(\mathbb{D}_{\frac{3}{4}}, \mathcal{M}_n)$, such that 

$$div (A \tilde{Q} \nabla \tilde{u} ) = \nabla A \tilde{Q}\nabla \tilde{u} +  A \nabla^\perp B (\tilde{Q} \nabla \tilde{u})= \nabla^\bot C \nabla \tilde{u} .$$
Hence $A$ and $C$ must satisfy
\begin{equation}
\left\{
\begin{array}{ll}
\Delta A &= \nabla^\bot C \nabla (\tilde{Q}^{-1})-\nabla A \nabla^\bot B \\
\Delta C &= \nabla A \nabla^\bot \tilde{Q} -div(A\nabla B \tilde{Q})
\end{array}
\right.
\end{equation}
This is exactly the same system as in (V.36) of \cite{RIV2012}. Hence there exists $A\in W^{1,2}(\mathbb{D}_{\frac{3}{4}}, Gl(m))$ and $C \in W^{1,2}(\mathbb{D}_{\frac{3}{4}}, \mathcal{M}_m)$ such that
\begin{equation}
\left\{
\begin{array}{ll}
\Delta A &= \nabla^\bot C \nabla (\tilde{Q}^{-1})-\nabla A \nabla^\bot B  \hbox{ on } \mathbb{D}_{\frac{3}{4}} \\
\partial_\nu A &=0 \hbox{ on } \partial \mathbb{D}_{\frac{3}{4}} \hbox{ and } \int_{\mathbb{D}_{\frac{3}{4}}} A =\frac{9\pi}{16} Id\\
\Delta C &= \nabla A \nabla^\bot \widetilde{Q} +div(A\nabla B \widetilde{Q})
\hbox{ on } \mathbb{D}_{\frac{3}{4}} \\
C&=0  \hbox{ on } \partial \mathbb{D}_{\frac{3}{4}}
\end{array}
\right.
\end{equation}
with
$$ \int_{\mathbb{D}_{\frac{3}{4}}} \vert \nabla A \vert^2 \, dx + \Vert dist(A, SO(m))\Vert_\infty \leq C \int_{\mathbb{D}_\frac{3}{4}} \vert \Omega \vert^2 \, dx$$
and

$$\int_{\mathbb{D}_{\frac{3}{4}}} \vert \nabla C \vert^2 \, dx \leq C \int_{\mathbb{D}_\frac{3}{4}} \vert \Omega  \vert^2 \, dx .$$

Finally setting, $\widetilde{A}=A\widetilde{Q}$, we have 
$$\Delta (\widetilde{A} \nabla \widetilde{u}) =\nabla^\bot C \nabla \widetilde{u} .$$

Then we are reduced to some classical Wente-type equation and the result follow directly from theorem V.3 of \cite{RIV2012}. We would like to mention that we have to pay attention, that the Wente inequality is in general false with the Neumann data, but here since $ C \in W^{1,2}_0$ it works, see\cite{DALP}.\hfill$\diamondsuit$\\

\bibliographystyle{plain}
\bibliography{biblio}

\end{document}